\documentclass[11pt]{article}
\usepackage{amsmath}
\usepackage{amsfonts,amssymb,amstext}
\usepackage{amsmath}
\usepackage{amstext}
\usepackage{amsthm}
\usepackage{graphicx}
\usepackage{setspace}
\usepackage{xcolor}
\usepackage{booktabs}

\makeatletter

\newcommand{\beq}{\begin{equation}}
\newcommand{\eeq}{\end{equation}}
\newcommand{\bea}{\begin{eqnarray}}
\newcommand{\eea}{\end{eqnarray}}
\newcommand{\beaa}{\begin{eqnarray*}}
\newcommand{\eeaa}{\end{eqnarray*}}

\newcommand{\n}{\noindent}
\newcommand{\q}{\quad}
\newcommand{\qq}{\qquad}

\newcommand{\E}{\rm I\!E}

\newcommand{\g}{\gamma}
\newcommand{\I}{\varphi}
\newcommand{\G}{\Gamma}
\newcommand{\de}{\delta}
\newcommand{\De}{\Delta}
\newcommand{\La}{\Lambda}
\newcommand{\al}{\alpha}
\newcommand{\la}{\lambda}
\newcommand{\f}{\infty}
\newcommand{\vs}{\varepsilon}
\newcommand{\cd}{\cdot}
\newcommand{\si}{\sigma}

\newcommand{\be}{\beta}

\newcommand{\Om}{\Omega}
\newcommand{\om}{\omega}

\newcommand{\st}{\subset}

\newcommand{\inl}{\int_{-\pi}^{\pi}}

\newcommand {\ol} {\overline}

\newcommand {\ra} {\rightarrow}

\newcommand {\s} {\section}
\newcommand {\sn} {\subsection}
\newcommand {\ssn} {\subsubsection}

\newtheorem{thm}{Theorem}[section]

\newtheorem{pp}{Proposition}[section]
\newtheorem{cor}{Corollary}[section]

\newtheorem{exa}{Example}[section]
\newtheorem{rem}{Remark}[section]

\newtheorem{den}{Definition}[section]

\numberwithin{equation}{section}

\usepackage[active]{srcltx}

\begin{document}

\title{Asymptotic behavior of the prediction error
for stationary sequences\thanks{The work is
dedicated to our teacher Academician, Professor  Il'dar Abdullovich 
Ibragimov on the occasion of his 90th birthday.}}

\author{Nikolay M. Babayan \thanks{ Russian-Armenian University, Yerevan, Armenia, e-mail: nmbabayan@gmail.com}
and Mamikon S. Ginovyan\thanks{Boston University, Boston, USA, e-mail: ginovyan@math.bu.edu}}

\date{}
\maketitle

\begin{abstract}
\noindent
One of the main problem in prediction theory of discrete-time second-order stationary
processes $X(t)$ is to describe the asymptotic behavior of the best linear mean
squared prediction error in predicting $X(0)$ given $ X(t),$ $-n\le t\le-1$,
as $n$ goes to infinity.
This behavior depends on the regularity (deterministic or nondeterministic)
and on the dependence structure of the underlying observed process $X(t)$.
In this paper we consider this problem both for deterministic and
nondeterministic processes and survey some recent results.
We focus on the less investigated case - deterministic processes.
It turns out that for nondeterministic processes the asymptotic behavior
of the prediction error is determined by the dependence structure
of the observed process $X(t)$ and the differential properties of its
spectral density $f$, while for deterministic processes it is determined
by the geometric properties of the spectrum of $X(t)$ and singularities of its spectral density $f$.
\end{abstract}



\tableofcontents

\s{Introduction}
\label{Int}

\sn{The finite prediction problem}
Let $X(t),$ $t\in\mathbb{Z}: = \{0,\pm1,\ldots\}$, be a centered discrete-time
second-order stationary process. The process is assumed to have an absolutely
continuous spectrum with spectral density function $f(\la),$
$\la\in [-\pi, \pi].$ The {\it 'finite' linear prediction problem}
is as follows.

Suppose we observe a finite realization of the process $X(t)$:
$$\{X(t), \,\, -n\le t\le-1\}, \q n\in\mathbb{N}: = \{1,2, \ldots\}.$$
We want to make an one-step ahead prediction, that is, to predict the
unobserved random variable $X(0)$,  using the {\it linear predictor}

$$Y=\sum_{k=1}^{n}c_kX(-k).$$

The coefficients
$c_k$, $k=1,2,\ldots,n$, are chosen so as to minimize
{\it the mean-squared error}: $\E\left|X(0) - Y\right|^2,$
where $\E[\cd]$ stands for the expectation operator.
If such minimizing constants $\widehat c_k:=\widehat c_{k,n}$
can be found, then the random variable
$$\widehat X_n(0):=\sum_{k=1}^{n}\widehat c_kX(-k)$$
is called {\it the best linear one-step ahead predictor} of $X(0)$
based on the observed finite past:\\ $X(-n), \ldots, X(-1)$.
The minimum mean-squared error:
$$\si_{n}^2(f): =\E\left|X(0) - \widehat X_n(0)\right|^2\geq0
$$
is called {\it the best linear one-step ahead prediction error} of $X(t)$
based on the past of length $n$.

One of the main problem in prediction theory of second-order stationary
processes, called {\it the 'direct' prediction problem} is to describe the
asymptotic behavior of the prediction error $\sigma_n^2(f)$ as $n\to\f$.
This behavior depends on the regularity nature (deterministic or nondeterministic)
of the observed process $X(t)$.

Observe that $\sigma_{n+1}^2(f)\leq \sigma_n^2(f)$, $n\in\mathbb{N}$,
and hence the limit of $\sigma_n^2(f)$ as $n\to\f$ exists.
Denote by $\sigma^2(f): = \sigma_\infty^2(f)$ the prediction error
of $X(0)$ by the entire infinite past: $\{X(t)$, $t\le-1\}$.

From the prediction point of view it is natural to distinguish the class of
processes for which we have {\it error-free prediction} by the entire infinite past, that is,
$\si^2(f)=0$. Such processes are called {\it deterministic} or {\it singular}.
Processes for which $\si^2(f)>0$ are called {\it nondeterministic}.

\vskip1mm
\n
{\it Note.} The term "deterministic" here is not used in the usual sense of
absence of randomness.
Instead determinism of a process means that there is an extremely strong dependence
between the successive random variables forming the process, yielding  error-free prediction
when using the entire infinite past (for more about this term see
Section \ref{dnd}, and also Bingham \cite{Bin1},
and Grenander and Szeg\H{o} \cite{GS}, p.176.)

Define the {\it 'relative' prediction error}
$$\de_n(f): = \si^2_n(f) - \si^2(f),$$
and observe that
$\delta_n(f)$ is non-negative and tends to zero as $n\to\infty$.
But what about the speed of convergence of $\delta_n(f)$ to zero as $n\to\infty$?
The paper deals with this question.
Specifically, the prediction problem we are interested in is
{\it to describe the rate of decrease of $\delta_n(f)$ to zero as}
$n \to \infty,$ depending on the regularity nature
of the observed process $X(t)$.

We consider the problem both for deterministic and
nondeterministic processes and survey some recent results.
We focus on the less investigated case - deterministic processes.
It turns out that for nondeterministic processes the asymptotic behavior
of the prediction error is determined by the dependence structure
of the observed process $X(t)$ and the differential properties of its
spectral density $f$, while for deterministic processes it is determined
by the geometric properties of the spectrum of $X(t)$ and singularities of its spectral density $f$.

\sn{A brief history}
The prediction problem stated above goes back to classical works of
A. N. Kolmogorov \cite{Kol1,Kol2}, G. Szeg\H{o} \cite{Sz,Sz1}
and N. Wiener \cite{Wr1}.
It was then considered by many authors for different classes of nondeterministic
processes (see, e.g., Baxter \cite{Bax}, Devinatz \cite{Dev1},
Geronimus \cite{Ger-3,Ger-1}, Golinski \cite{Gol-1},
Golinski and Ibragimov \cite{GI}, Grenander and Rosenblatt \cite{GR1, GR2}, Grenander
and Szeg\H{o} \cite{GS}, Helson and Szeg\H{o} \cite{HS}, Hirshman \cite{Hir},
Ibragimov \cite{I-2,I-4}, Ibragimov and Rozanov \cite{IR}, Ibragimov and Solev \cite{IS}, Inoue \cite{In2},
Pourahmadi \cite{Po}, Rozanov \cite{R}, and reference therein).
More references can be found in the survey papers Bingham \cite{Bin1} and
Ginovyan \cite{G-4}.
In Section \ref{ndp} of the paper we state some important known results for
nondeterministic processes.

We focus in this paper on deterministic processes, that is, on the
class of processes for which  $\si^2(f) =0$.
This case is not only of theoretical interest, but is also important from the point
of view of applications.
For example, as pointed out by Rosenblatt \cite{Ros} (see also Pierson \cite{Pl}),
situations of this type arise in Neumann's theoretical model of
storm-generated ocean waves.
Such models are also of interest
in meteorology (see, e.g., Fortus \cite{For}).

Only few works are devoted to the study of the speed of convergence of
$\de_n(f)=\si^2_n(f)$ to zero as $n\to\infty$, that is, the asymptotic behavior of
the prediction error for deterministic processes.
One needs to go back to the classical work of M. Rosenblatt \cite{Ros}.
Using the technique of orthogonal polynomials on the unit circle (OPUC),
M. Rosenblatt investigated the asymptotic behavior of the prediction
error $\si^2_n(f)$ for deterministic processes in the following two cases:
\begin{itemize}
\item[(a)]
the spectral density $f(\la)$ is continuous and positive on a
segment of $[-\pi,\pi]$ and is zero elsewhere,

\item[(b)]
the spectral density $f(\la)$ has a very high order of contact with zero
at points $\la=0,\pm\pi$, and is strictly positive otherwise.
\end{itemize}

Later the problems (a) and (b) were studied by  Babayan  \cite{Bb-1,Bb-2}, Babayan and Ginovyan \cite{BG-1, BG, BG-2}, Babayan et al. \cite{BGT}
(see also Davisson \cite{Dav1} and Fortus \cite{For}), where some
generalizations and
extensions of Rosenblatt's results have been obtained.

\sn{Notation and conventions.}

Throughout the paper we will use the following notation and conventions.\\
The standard symbols $\mathbb{N}$, $\mathbb{Z}$, $\mathbb{R}$ and $\mathbb{C}$
denote the sets of natural, integer, real and complex numbers, respectively.
Also, we denote
$\mathbb{Z}_+: = \{0,1,2,\ldots\}$, $\Lambda: = [-\pi, \pi],$
$\mathbb{D}:=\{z\in\mathbb{C}: |z|<1\}$, $\mathbb{T}:=\{z\in \mathbb{C}: \, |z|=1\}$.
For a point $\la_0\in \Lambda$ and a number $\de>0$ by $O_\de(\la_0)$ we
denote a $\de$-neighborhood of $\la_0$, that is,
$O_\de(\la_0):=\{\la\in \Lambda: \ |\la-\la_0|<\de\}$.
By $L^p(\mu):=L^p(\mathbb{T},\mu)$ ($p\geq $1) we denote the weighted
Lebesgue space with respect to the measure $\mu$, and by $\|\cdot\|_{p,\mu}$ we denote the norm in $L^p(\mu)$.
In the special case where $\mu$ is the Lebesgue measure,
we will use the notation $L^p$ and $||\cdot||_{p}$, respectively.
For a function $h\geq0$ by $G(h)$ we denote the geometric mean of $h$.
For two functions $f(\la)\geq0$ and $g(\la)\geq0$ 
we will write $f(\lambda){\sim}g(\lambda)$ as ${\lambda\to\lambda_0}$
if $\lim_{\lambda\to\lambda_0}{f(\lambda)}/{g(\lambda)}=1$;
$f(\lambda){\simeq}g(\lambda)$ as ${\lambda\to\lambda_0}$ if
$\lim_{\lambda\to\lambda_0}{f(\lambda)}/{g(\lambda)}=c>0$,
and $f(\lambda){\asymp}g(\lambda)$ if there are constants $c_1, c_2$
($0<c_1\leq c_2<\f$) such that $0<c_1\leq f(\lambda)/g(\lambda)
\leq c_2<\f$ for all $\la\in\Lambda$.
We will use similar notation for sequences: for two sequences
$\{a_n\geq0, n\in\mathbb{N}\}$ and $\{b_n>0, n\in\mathbb{N}\}$,
we will write $a_n\sim b_n$ if $\lim_{n\to\f}{a_n}/{b_n}=1$,
$a_n{\simeq} b_n$ if $\lim_{n\to\f}{a_n}/{b_n}=c>0$,
$a_n{\asymp} b_n$ if $c_1\leq {a_n}/{b_n}\leq c_2$ for all $\la\in\mathbb{N}$,
$a_n=O(b_n)$ if ${a_n}/{b_n}$ is bounded,
and $a_n=o(b_n)$ if ${a_n}/{b_n}\to0$ as $n\to\f$.
For a set $E$ by $\ol E$ we denote the closure of $E$.
The letters $C$, $c$, $M$ and $m$ with or without indices are used
to denote positive constants, the values of which can vary from line to line.

We will use the abbreviations: OPUC for 'orthogonal polynomials
on the unit circle', PACF for 'partial autocorrelation function',
and 'a.e.' for 'almost everywhere' (with respect to the Lebesgue measure).
We will assume that all the relevant objects are defined in terms
of Lebesgue integrals,
and so are invariant under change of the integrand on a null set.

\sn{The structure of the paper}

The paper is structured as follows.
In Section \ref{model} we describe the model of interest - a stationary process,
and recall some key notions and results from the theory of stationary processes.
In Section \ref{fpp} we present formulas for
the finite prediction error $\si_{n}^2(f)$, and state some preliminary results.
In Section \ref{ndp} we state some well known results on asymptotic behavior
of the prediction error for nondeterministic processes.
Asymptotic behavior of the finite prediction error $\si_{n}^2(f)$ for deterministic processes is
discussed in Section \ref{d}. Here we state extensions of Rosenblatt's
and Davisson's results, and discuss a number of examples.
In Section \ref{Ap} we analyze the relationship between the rate of
convergence to zero of the prediction error $\si_{n}^2(f)$ and the
minimal eigenvalue of a truncated Toeplitz matrix generated by the
spectral density $f$.
In Section \ref{methods} we briefly discuss the
tools, used to prove
the theorems stated in the paper.

\section{The model. Key notions and some basic results}
\label{model}
In this section we introduce the model of interest - a second-order
stationary process, and recall some key notions and results from the
theory of stationary processes.

\sn{Second-order (wide-sense) stationary processes}
Let $\{X(t), \ t\in \mathbb{Z}\}$ be a centered real-valued second-order
(wide-sense) stationary process defined on a probability space
$(\Omega, \mathcal{F}, P)$ with covariance function $r(t)$, that is,
\[
{\E}[X(t)]=0, \q r(t)={\E}[X(t+s)X(s)], \q s, t\in\mathbb{Z},
\]
where ${\E}[\cdot]$ stands for the expectation operator with respect
to the measure $P$.

By the Herglotz theorem (see, e.g., Brockwell and Davis \cite{BD},
p. 117-118),
there is a finite measure $\mu$ on $\Lambda$ such that
the covariance function $r(t)$
admits the following {\it spectral representation}:
\beq
\label{i1}
r(t)=\inl e^{-it\la}d\mu(\la), \q t\in \mathbb{Z}.
\eeq
\n
The measure $\mu$ in (\ref{i1}) is called the {\it spectral measure} of the
process $X(t)$.
If $\mu$ is absolutely continuous (with respect to the Lebesgue measure),
then the function $f(\la):=d\mu(\la)/d\la$ is called the {\it spectral density} of $X(t)$.
We assume that $X(t)$ is a {\it non-degenerate} process, that is,
${\rm Var}[X(0)]:={\E}|X(0)|^2=r(0)>0$ and, without loss of generality,
we may take $r(0)=1$.
Also, to avoid the trivial cases, we assume that the spectral measure
$\mu$ is {\it non-trivial}, that is, $\mu$ has infinite support.

Notice that if the spectral density $f(\la)$ exists, then
$f(\la)\geq 0$, $f(\la)\in L^1(\Lambda)$, and \eqref{i1} becomes
\begin{equation}
\label{mo1}
r(t)=\inl e^{-it\la}f(\la)d\la, \q t\in \mathbb{Z}.
\end{equation}
Thus, the covariance function $r(t)$ and the spectral
function $F(\la)$ (resp. the spectral density function $f(\la)$) are equivalent
specifications of the second order properties for a stationary process
$\{X(t), \ t\in\mathbb{Z}\}$.

\begin{rem}
{\rm The parametrization of the unit circle $\mathbb{T}$ by the formula
$z=e^{i\la}$ establishes a bijection between $\mathbb{T}$ and the
interval $[-\pi,\pi)$. By means of this bijection the measure $\mu$ on
$\Lambda$ generates the corresponding measure on the unit circle $\mathbb{T}$,
which we also denote by $\mu$. Thus, depending
on the context, the measure $\mu$ will be supported either on $\Lambda$ or
on $\mathbb{T}$.
We use the standard Lebesgue decomposition of the measure $\mu$:
\beq
\label{di1}
d\mu(\la)=d\mu_{a}(\la) +d\mu_s(\la) =f(\la)d\la +d\mu_s(\la),
\eeq
where $\mu_a$ is the absolutely continuous part of $\mu$
(with respect to the Lebesgue measure)
and $\mu_s$ is the singular part of $\mu$, which is the sum
of the discrete and continuous singular components of $\mu$.}
\end{rem}

By the well-known Cram\'er theorem (see, e.g., Cram\'er and Leadbetter \cite{CL}) 
for any stationary process $\{X(t), \, t\in \mathbb{Z}\}$ with spectral
measure $\mu$ there exists an orthogonal stochastic measure $Z=Z(B)$,
$B\in\mathfrak{B}(\Lambda)$, such that for every $t\in \mathbb{Z}$
the process $X(t)$ admits the following {\sl spectral representation}:
\begin{equation}
\label{i18}
X(t)=\int_\Lambda e^{-it\la}dZ(\la), \q  t\in\mathbb{Z}.
\end{equation}
Moreover, ${\E}\left[|Z(B)|^2\right]=\mu(B)$ for every $B\in\mathfrak{B}(\Lambda)$.
Here $\mathfrak{B}(\Lambda)$ stands for the Borel $\si$-algebra of the sets of $\Lambda$.
For definition and properties of orthogonal stochastic measures and
stochastic integral in (\ref{i18}) we refer, e.g., Cram\'er and Leadbetter \cite{CL}, Ibragimov and Linnik \cite{IL}, and Shiryaev \cite{Sh}.

\sn{Linear processes. Existence of spectral density functions}

We will consider here stationary processes possessing spectral density
functions. For the following results we refer to
Cram\'er and Leadbetter \cite{CL}, Doob \cite{Doob}, Ibragimov and Linnik \cite{IL}.
\begin{thm} The following assertions hold.
\begin{itemize}
\item[(a)]
The spectral function $F(\la)$ of a stationary process $\{X(t), \,t \in \mathbb{Z}\}$
is absolutely continuous (with respect to the Lebesgue measure), that is, $F(\la)=\int_{-\pi}^\la f(x)dx$,
if and only if it can be represented as an infinite moving average:
\begin{equation}
\label{dlp}
X(t) = \sum_{k=-\f}^{\f}a(t-k)\xi(k), \quad
         \sum_{k=-\f}^{\f}|a(k)|^2 < \f,
\end{equation}
where $\{\xi(k), k\in \mathbb{Z}\}\sim$ WN(0,1) is a standard white-noise, that is, a sequence of orthonormal random variables.
\item[(b)]
The covariance function $r(t)$ and the spectral density $f(\la)$ of $X(t)$ are given by formulas:
\begin{equation}
\label{dcv}
r(t)= {\E} X(t)X(0)=\sum_{k=-\f}^{\f}a(t+k) a(k),
\end{equation}
and
\beq
\label{dsd}
f(\la) =
\frac{1}{2\pi}\left|\sum_{k=-\f}^{\f}a(k)e^{-ik\la}\right|^2
= \frac{1}{2\pi} |\widehat a(\la)|^2,
\quad \la\in\Lambda.
\eeq
\item[(c)]
In the case where $\{\xi(k), k\in \mathbb{Z}\}$ is a sequence of Gaussian random variables,
the process $\{X(t), \,t \in \mathbb{Z}\}$ is Gaussian.
\end{itemize}
\end{thm}

\sn{Dependence (memory) structure of the model}
\label{mem}

Depending on the memory (dependence) structure, we will distinguish
the following types of stationary models:

(a) short memory (or short-range dependent),

(b) long memory (or long-range dependent),

(c) intermediate memory (or anti-persistent).

\noindent
The memory structure of a stationary process is essentially a
measure of the dependence between all the variables in the process,
considering the effect of all correlations simultaneously.
Traditionally memory structure has been defined in the time domain
in terms of decay rates of the autocorrelations, or in the
frequency domain in terms of rates of explosion of low frequency spectra
(see, e.g., Beran et al. \cite{BFGK}, and references therein).

It is convenient to characterize the memory structure in terms
of the spectral density function.

\ssn{Short memory models}
A stationary process $\{X(t), \,t \in \mathbb{Z}\}$ with spectral
density function $f(\la)$ is said to be a {\it short memory} process
if the spectral density $f(\lambda)$ is bounded away from zero and infinity,
that is, there are constants $C_1$ and $C_2$ such that
\begin{equation*}
0< C_1 \le f(\la) \le C_2 <\f.
\end{equation*}
A typical short memory model example is the stationary
Autoregressive Moving Average (ARMA)$(p,q)$ process $X(t)$
defined to be a stationary solution of the difference equation:
\begin{equation*}
\psi_p(B)X(t)=\theta_q(B)\varepsilon(t), \q t\in\mathbb{Z},
\end{equation*}
where $\psi_p$ and $\theta_q$ are polynomials respectively
of degrees $p$ and $q$ having no zeros on the unit circle $\mathbb{T}$,
$B$ is the backshift operator defined by $BX(t)=X(t-1)$,
and $\{\vs(t), t\in\mathbb{Z}\}$ is a WN(0,$\si^2$) white noise,
that is, a sequence of zero-mean, uncorrelated random variables
with variance $\si^2$.
The covariance $r(k)$ of (ARMA)$(p,q)$ process is
exponentially bounded:
\[
|r(t)|\le Cr^{-t}, \q t=1,2,\ldots;\q0<C<\f; \,\, 0<r<1,
\]
and the spectral density $f(\la)$ is a rational function
(see, e.g., Brockwell and Davis \cite{BD}, Sec. 3.1):
%
\begin{equation}
\label{arma}
f(\la) = \frac{\si^2}{2\pi}\cd\frac
{|\theta_q(e^{-i\la})|^2}{|\psi_p(e^{-i\la})|^2}.
\end{equation}

\ssn{Long-memory and anti-persistent models}

A {\it long-memory} model is defined to be a
stationary process with {\it unbounded} spectral density,
and an {\it anti-persistent} model -- a stationary
process with {\it vanishing} (at some fixed points) spectral density
(see, e.g., Beran et al. \cite{BFGK}, Brockwell and Davis \cite{BD},
and references therein).

A typical model example that displays long-memory and intermediate
memory (anti-persistent) is the Autoregressive Fractionally
Integrated Moving Average (ARFIMA)$(p,d,q)$ process $X(t)$ defined
to be a stationary solution of the difference equation
(see, e.g., Brockwell and Davis \cite{BD}, Section 13.2):
\begin{equation*}
\psi_p(B)(1-B)^dX(t)=\theta_q(B)\varepsilon(t), \q d<1/2,
\end{equation*}
where $B$ is the backshift operator, $\varepsilon(t)$ is a
WN(0,$\si^2$) white noise, and $\psi_p$ and $\theta_q$ are
polynomials of
degrees $p$ and $q$, respectively.
The spectral density $f_X(\la)$ of $X(t)$ is given by
\begin{equation}
\label{AA}
f_X(\la)=|1-e^{-i\la}|^{-2d}f(\la)=(2\sin(\la/2))^{-2d}f(\la), \q d<1/2,
\end{equation}
where $f(\la)$ is the spectral density of an ARMA$(p,q)$ process,
given by (\ref{arma}). The condition $d<1/2$ ensures that
$\int_{-\pi}^\pi f(\la)d\la<\f$,
implying that the process $X(t)$ is well defined because
${\E}|X(t)|^2=\int_{-\pi}^\pi f(\la)d\la.$

\noindent
Observe that for $ 0<d<1/2$ the model $X(t)$ specified by the spectral
density (\ref{AA}) displays long-memory. In this case we have
$f(\la)\thicksim c\, |\la|^{-2d}$ as $\la\to0$, that is, $f(\la)$ blows up
at $\la=0$ like a power function, which is the typical behavior of
a long memory model. For $d<0$, the model $X(t)$ displays
intermediate-memory, and in this case, the spectral density in \eqref{AA}
vanishes at $\la=0$. For $d=0$ the model $X(t)$ displays short-memory.
For $d \ge1/2$ the function $f_X(\la)$ in (\ref{AA}) is not integrable,
and thus it cannot represent a spectral density of a stationary process.

\sn {Deterministic and nondeterministic processes.}
\label{dnd}
In this section we state Kolmogorov's isometric isomorphism theorem
and the infinite prediction problem. We give time-domain
(Wold's theorem) and
frequency-domain (Kolmogorov-Szeg\H{o}'s theorem) characterizations of deterministic and nondeterministic processes.

\ssn{Kolmogorov's isometric isomorphism theorem}
Given a probability space $(\Om, \mathcal{F}, P)$, define the $L^2$-space of real-valued random variables $\xi= \xi(\om)$ with
$\E[\xi]=0$:
\begin{equation}
\label{i19}
L^2(\Om):=\left\{\xi: \, ||\xi||^2:=
\int_\Om |\xi(\om)|^2dP(\om)<\f\right\}.
\end{equation}
Then $L^2(\Om)$ becomes a Hilbert space with the following
inner product: for $\xi, \eta\in L^2(\Om)$
\begin{equation}
\label{i20}
(\xi,\eta)= \E[\xi \eta]=\int_\Om \xi(\om){\eta(\om)}dP(\om).
\end{equation}
\n
For $a, b\in \mathbb{Z}$,  $-\f\le a\le b\le \f,$
we define  the space $H_a^b(X)$ to be the closed linear
subspace of the space $L^2(\Om)$ spanned by the random
variables $X(t)=X(t,\om)$, $t\in [a,b]$, $\om\in\Omega$:
\begin{equation}
\label{i21}
 H_{a}^b(X): = \ol{sp}\{X(t), \,\, a \le t\le b\}_{L^2(\Om)}.
\end{equation}

\n Observe that the space $H_{a}^b(X)$ consists of
all finite linear combinations of the form $\sum_{k=a}^b c_k X(k)$,
as well as, their $L^2(\Om)$-limits.\\
The space $H(X):=H_{-\f}^\f(X)$ is called the {\it Hilbert
space generated by the process} $X(t)$, or the {\it time-domain} of $X(t)$.

Let $\mu$ be the spectral measure of the process $\{X(t), \, t\in \mathbb{Z}\}$.
Consider the weighted $L^2$-space $L^2(\mu):=L^2(\mu,\Lambda)$ of complex-valued functions
$\I(\la), \, \la\in\Lambda$, defined by
\beq
\label{i22}
L^2(\mu):=\left\{\I(\la): \,||\I||^2_\mu:= \int_\Lambda|\I(\la)|^2d\mu(\la)<\f\right\}.
\eeq
Then $L^2(\mu)$ becomes a Hilbert space with the following
inner product: for $\I, \psi\in L^2(\mu)$
\begin{equation}
\label{i23}
(\I, \psi)_\mu= \int_\Lambda\I(\la)\ol{\psi}(\la)d\mu(\la).
\end{equation}
The Hilbert space $L^2(\mu,\Lambda)$ is called the {\it frequency-domain} of the process $X(t)$.

\begin{thm}[Kolmogorov's isometric isomorphism theorem]
\label{Kolm}
For any stationary process $X(t)$, $t\in \mathbb{Z}$, with spectral
measure $\mu$ there exists a unique isometric isomorphism $V$ between
the time-domain $H(X)$ and the  frequency-domain $L^2(\mu)$,
such that $V[X(t)]=e^{it\la}$ for any $t\in \mathbb{Z}.$
\end{thm}
Thus, in view of Theorem \ref{Kolm}, any linear problem in
the time-domain $H(X)$
can be translated into one in the frequency-domain
$L^2(\mu)$, and vice versa. This fact allows to
study stationary processes using analytic methods.

\ssn{The infinite prediction problem}
\label{ipp}

\n Observe first that since by assumption $X(t)$ is a non-degenerate
process, the time-domain $H(X)$ of $X(t)$ is non-trivial, that is,
$H(X)$ contains elements different from zero.
\begin{den}
{\rm The space $H_{-T}^t(X)$ is called the
{\it finite history}, or {\it past of length $T$ and present}
of the process $X(u)$ up to time $t$.
\n The space $H_t(X):=H_{-\f}^t(X)$ is called the
{\it entire history}, or {\it infinite past and present}
of the process $X(u)$ up to time $t$.
\n The space
\beq
\label{w1}
H_{-\f}(X): = \cap _t H_{-\f}^t (X)
\eeq
is called the {\it remote past} of the process $X(u)$.}
\end{den}
\n
It is clear that
\beq
\label{w2}
H_{-\f}(X) \st \cdots \st H_{-\f}^t(X) \st H_{-\f}^{t+\tau}(X)
\st \cdots \st H(X), \q \tau\in\mathbb{N}.
\eeq

The Hilbert space setting provides a natural framework for
stating and solving the problem of predicting future values
of the process $X(u)$ from the observed past values.
Assume that a realization of the process $X(u)$ for times $u\le t$
is observed and we want to predict the value $X(t+\tau)$ for some
$\tau\ge1$ from the observed values. Since we will never
know what particular realization is being observed, it is
reasonable to consider as a predictor $\widehat X(t,\tau)$
for $X(t+\tau)$ a function of the observed values,
$g(\{X(u),u\le t\})$, which is good "on the average".
So, as an optimality criterion for our predictor we take
the $L^2$-distance, that is, the mean squared error, and
consider only the linear predictors. With these restrictions,
the infinite linear prediction problem can be stated as follows.

\vskip2mm
\n {\bf The infinite linear prediction problem.}
Given a 'parameter' of the process $X(u)$  (e.g., the covariance
function $r(t)$ or the spectral function $F(\la)$),
the entire history $H_{-\f}^t(X)$ of $X(u)$, and a number $\tau\in\mathbb{N}$,
find a random variable $\widehat X(t,\tau)$ such that
\begin{itemize}
\item[a)]
$\widehat X(t,\tau)$ is {\it linear}, that is,
$\widehat X(t,\tau)\in H_{-\f}^t(X)$,
\item[b)]
$\widehat X(t,\tau)$ is {\em mean-square optimal (best)}
among all elements $Y\in H_{-\f}^t(X)$, that is,
$\widehat X(t,\tau)$ minimizes the mean-squared error
$||X(t+\tau)-Y||^2_{L^2(\Om)}:$
\beq
\label{w3}
||X(t+\tau)-\widehat X(t,\tau)||^2_{L^2(\Om)}=
\min_{Y\in H_{-\f}^t(X)}||X(t+\tau)-Y||^2_{L^2(\Om)}.
\eeq
\end{itemize}

\n The solution - the random variable $\widehat X(t,\tau)$
satisfying a) and b), is called the
{\it best linear $\tau$-step ahead predictor} for an element
$X(t+\tau)\in H(X)$. The quantity
\beq
\label{w4}
\si^2(\tau):= ||X(t+\tau)-\widehat X(t,\tau)||^2_{L^2(\Om)}
= ||X(t+\tau)||^2_{L^2(\Om)}-||\widehat X(t,\tau)||^2_{L^2(\Om)},
\eeq
which is independent of $t$, is called the {\it prediction error (variance)}.

\vskip2mm

The advantage of the Hilbert space setting now becomes apparent.
Namely, by the {\it projection theorem} in Hilbert spaces
(see, e.g., Pourahmadi \cite{Po}, p. 312),
such a predictor $\widehat X(t,\tau)$ exists, is unique, and is given by
\beq
\label{w5}
\widehat X(t,\tau)= P_t X(t+\tau),
\eeq
where $P_t:=P_{-\f}^t$ is the orthogonal projection
operator in $H(X)$ onto $H_{-\f}^t(X)$.
\begin{rem}
\label{r311}
{\rm The reason for restricting attention to linear predictors
is that the best linear predictor $\widehat X(t,\tau)$,
in this case, depends only on knowledge of the covariance
function $r(t)$ or the spectral function $F(\la)$.
The prediction problem becomes much more difficult when
nonlinear predictors are allowed (see, e.g., Hannan \cite{Han},
Koopmans \cite{Kp}).}
\end{rem}

\ssn {Deterministic (singular) and nondeterministic processes. Characterizations}

From prediction point of view it is natural to distinguish
the class of processes for which we have {\it error-free prediction},
that is, $\si(\tau)=0$  for all $\tau\ge1$, or equivalently,
$\widehat X(t,\tau)=X(t+\tau)$ for all $t\in \mathbb{Z}$ and
$\tau\ge1$.
In this case, the prediction is called {\it perfect}.
It is clear that a process $X(t)$ possessing perfect prediction
represents a singular case of {\it extremely strong dependence}
between the random variables forming the process.
Such a process $X(t)$ is called {\it deterministic} or {\it singular}.
From the physical point of view, singular processes are exceptional.
From application point of view, of considerable interest is the class of processes for which we have $\si(\tau)>0$ for all $\tau\ge1$.
In this case the prediction is called {\it imperfect}, and the
process $X(t)$ is called {\it nondeterministic}.

Observe that the time-domain $H(X)$ of any non-degenerate
stationary process $\{X(t)$, $t\in\mathbb{Z}\}$ can be
represented as the orthogonal sum $H(X)= H_1(X)\oplus H_{-\f}(X),$
where $H_{-\f}(X)$ is the remote past of $X(t)$ defined by (\ref{w1}),
and $H_1(X)$ is the orthogonal complement of $H_{-\f}(X)$.
So, we can give the following geometric definition of the deterministic
(singular), nondeterministic and purely nondeterministic (regular) processes.
\begin{den}{\rm A stationary process $\{X(t)$, $t\in\mathbb{Z}\}$ is called}
\begin{itemize}
\item
deterministic or singular {\rm if $H_{-\f}(X)= H(X),$ that is,
$H_{-\f}^t(X)=H_{-\f}^s(X)$ for all $t, s \in\mathbb{Z}$},
\item
nondeterministic {\rm if $H_{-\f}(X)$ is a proper subspace
of $H(X)$, that is, $H_{-\f}(X)\subset H(X)$},
\item
purely nondeterministic (PND) or regular {\rm if $H_{-\f}(X)=\{0\}$,
that is, the remote past $H_{-\f}(X)$ of $X(t)$ is the trivial subspace,
consisting of the singleton zero.}
\end{itemize}
\end{den}

\n The next theorem contains a characterization of deterministic
and purely nondeterministic processes in terms of prediction error.
\begin{thm}
\label{t30}
A stationary process $\{X(t)$, $t\in\mathbb{Z}\}$ is
\begin{itemize}
\item[(a)]
deterministic if and only if $\si(\tau_0) = 0$ for some $\tau_0\ge1$,
$\tau_0\in\mathbb{N}$ (then $\si(\tau) = 0$ for all $\tau\in\mathbb{N}$).
\item[(b)]
purely nondeterministic if and only if
$\lim_{\tau\to\f}\si^2(\tau) = E|X(t)|^2 = r(0)(=1).$
\end{itemize}
\end{thm}

\begin{rem}{\rm Every purely nondeterministic process $X(t)$
is nondeterministic, but the converse is generally not true.
An example of such process provides the process $X(t)=\vs(t)+\xi$,
where $\{\vs(t), \,t\in \mathbb{Z}\}$ is a ${\rm WN}(0,\si^2)$
(a centered white noise with variance $\si^2$), and
$\xi$ is a random variable such that $Var(\xi)=\si^2_0$ and
$(\vs(t), \xi)=0$ for all $t\in \mathbb{Z}$
(see Pourahmadi \cite{Po}, p. 163).
Observe also that the only process $X(t)$ that is both
deterministic and purely nondeterministic is the
degenerate process. Assuming that $X(t)$ is a non-degenerate
process, we exclude this trivial case.}
\end{rem}

The next result, known as Wold's decomposition theorem, provides a key
step for solution of the infinite prediction problem in the time-domain
setting, and essentially says that any stationary process can be
represented in the form of a sum of two orthogonal stationary components,
one of which is perfectly predictable (singular component),
while for the other (regular component) an explicit formula
for the predictor can be obtained.

\begin{thm}[Wold's decomposition]
\label{th32}
Every centered non-degenerate discrete-time stationary
process $X(t)$ admits a unique decomposition: $$X(t) = X_S(t)+X_R(t),$$
where
\begin{itemize}
\item[(a)]
the processes $X_R(t)$ and $X_S(t)$ are stationary, centered, mutually
uncorrelated (orthogonal), and completely subordinated to $X(t)$, that is,
$H_{-\f}^t(X_R)\subseteq H_{-\f}^t(X)$ and $H_{-\f}^t(X_S)\subseteq H_{-\f}^t(X)$
for all $t\in\mathbb{Z}.$
\item[(b)]
the process $X_S(t)$ is deterministic (singular),
\item[(c)]
the process $X_R(t)$ is purely nondeterministic (regular) and has the
infinite moving-average representation:
\beq
\label{w8}
X_R(t) = \sum_{k=0}^{\f}a_k\vs_0(t-k), \q \sum_{k=0}^{\f}|a_k|^2 < \f,
\eeq
where $\vs_0(t)$ is an innovation of $X_R(t),$ that is, $\vs_0(t)$
is a standard white-noise process, such that
$H_{-\f}^t(X_R)=H_{-\f}^t(\vs_0)$ for all $t\in\mathbb{Z}$.
\end{itemize}
\end{thm}

The next result describes the
asymptotic behavior of the prediction error $\si_n^2(\mu)$ for a stationary
process $X(t)$ with spectral measure $\mu$ of the form \eqref{di1}
and gives spectral characterizations of deterministic, nondeterministic
and purely nondeterministic
processes (see, e.g., Grenander and Szeg\H{o} \cite{GS}, p. 44,
and Ibragimov and Rozanov \cite{IR}, p. 35-36).
\begin{thm}
\label{th34}
Let $X(t)$ be a non-degenerate stationary process with spectral
measure $\mu$ of the form \eqref{di1}. The following assertions hold.
\begin{itemize}
\item[(a)]{\rm (Kolmogorov-Szeg\H{o's} Theorem)}.
The following relations hold.
\bea
\label{c013}
\lim_{n\to\f}\si_n^2(\mu)=\lim_{n\to\f}\si_n^2(f)=\si^2(f)=2\pi G(f),
\eea
where $G(f)$ is the {\it geometric mean} of $f$, namely
\beq
\label{a2}
G(f): = \left \{
\begin{array}{ll}
\exp\left\{\frac1{2\pi}\inl\ln f(\la)\,d\la \right\} &
\mbox{if \, $\ln f \in {L}^1(\Lambda)$}\\
           0, & \mbox { otherwise,} \qq
           \end{array}
           \right.
\eeq

\item[(b)]
$
H_{-\f}^0(\mu_s)=H(\mu_s)
\Leftrightarrow  \sigma^2(\mu)=0 \Leftrightarrow
X(t) \,\,  is \,\, deterministic,
$
\item[(c)]
{\rm (Kolmogorov-Szeg\H{o}'s alternative)}.
Either
$$
H_{-\f}^0(\mu_a)=H(\mu_a) \Leftrightarrow \inl\ln f(\la)\,d\la = -\f
\Leftrightarrow  \sigma^2(f)=0 \Leftrightarrow
X(t) \,\,  is \,\, deterministic,
$$
or else
$$
H_{-\f}^0(\mu_a) \neq  H(\mu_a) \Leftrightarrow\inl\ln f(\la)\,d\la > -\f \Leftrightarrow
\sigma^2(f) >0 \Leftrightarrow  X(t) \,\, is \,\, nondeterministic.
$$
\item[(d)] The process $X(t)$ is regular (PND) if and only if
it is nondeterministic and $\mu_s\equiv0$.
\end{itemize}
\end{thm}

\begin{rem}
\label{r33}
{\rm The second equality in (\ref{c013}) was proved by Szeg\H{o} in 1920,
while the first equality was proved by Kolmogorov in 1941 (see, e.g.,
Hoffman \cite{Hof}, p. 49).}
It is remarkable that \eqref{c013} is independent of the singular part $\mu_s$.

\end{rem}

The condition $\ln f \in {L}^1(\Lambda)$ in \eqref{a2} is equivalent to
the {\it Szeg\H{o} condition}:
\beq
\label{S}
\inl\ln f(\la)\,d\la > -\f
\eeq
(this equivalence follows because $\ln f(\la)\le f(\la)$).
The Szeg\H{o} condition \eqref{S} is also called the {\it non-determinism condition}.

In this paper we consider the class of deterministic
processes with absolutely continuous spectra.

We will say that the spectral density $f(\la)$ has a
{\it very high order of contact with zero at a point} $\la_0$ if $f(\la)$
is positive everywhere except for the point $\la_0$, due to which
the Szeg\H{o} condition \eqref{S} is violated.
Observe that the Szeg\H{o} condition is related to the
character of the singularities (zeroes and poles) of the spectral density $f$,
and does not depend on the differential properties of $f$.
For example, for any  $a>0$, the function $\widehat f_a(\la)=\exp\{-|\la|^{-a}\}$
is infinitely differentiable.
In addition, for $a<1$ Szeg\H{o}'s condition is satisfied,
and hence the corresponding process $X(t)$ is nondeterministic,
while for $a\ge 1$ Szeg\H{o}'s condition is violated, and $X(t)$
is deterministic (see, e.g., Pourahmadi \cite{Po}, p.68,
Rakhmanov \cite{Ra4}).
Thus, according to the above definition, for $a\ge 1$ this function
has a very high order of contact with zero at the point $\la=0$.

\s{Formulas for the finite prediction error $\si_{n}^2(f)$ and some properties}
\label{fpp}

In this section, we provide various formulas for the prediction error $\sigma_n^2(f)$ in terms of orthogonal polynomials and their parameters
(Verblunsky's coefficients), Toeplitz determinants, state Szeg\H{o}'s,
Verblunsky's and Rakhmanov's theorems, and list a number of properties
of $\sigma_n^2(f)$.

\sn{Formulas for the prediction error $\si_{n}^2(f)$}
We present here formulas for the finite prediction error
$\si_{n}^2(f)$ and state some preliminary results, which will be used in the sequel.

Suppose we have observed the values $X(-n), \ldots, X(-1)$ of a centered,
real-valued stationary process $X(t)$ with spectral measure $\mu$
of the form \eqref{di1}.
The {\it one-step ahead linear prediction problem}
in predicting a random variable $X(0)$ based on the observed values
$X(-n), \ldots, X(-1)$ involves finding constants
$\widehat c_k:=\widehat c_{k,n}$, $k=1,2,\ldots,n$,
that minimize the one-step ahead prediction error: 
\bea
\label{OS1}
\si_{n}^2(\mu):
=\min_{\{c_k\}}\E\left|X(0) - \sum_{k=1}^{n}c_kX(-k)\right|^2
=\E\left|X(0) - \sum_{k=1}^{n}\widehat c_kX(-k)\right|^2.
\eea
Using Kolmogorov's isometric isomorphism $V:\, X(t)\leftrightarrow e^{it\la}$,
in view of (\ref{OS1}), for the prediction error $\sigma_n^2(\mu)$ we can write
\bea
\label{OS2}
\sigma_n^2(\mu)=
\min_{\{c_k\}}\left\|1 -\sum_{k=1}^{n}c_ke^{-ik\la}\right\|^2_{2,\mu} 
=\min_{\{q_n\in \mathcal{Q}_n\}}\left\|q_n\right\|^2_{2,\mu}, 
\eea
where $||\cdot||_{2,\mu}$ is the norm in $L^2(\mathbb{T},\mu)$, and
\beq
\label{Q_n}
\mathcal{Q}_n: =\left\{q_n:  q_n(z)=z^{n} + c_1z^{n-1}+ \cdots c_n\right\}
\eeq
is the class of monic polynomials (i.e. with $c_0=1$) of degree $n$.
Thus, the problem of finding $\sigma_n^2(\mu)$ becomes to the problem of finding
the solution of the minimum problem \eqref{OS2}-\eqref{Q_n}.

The polynomial $p_n(z):=p_n(z,\mu)$ which solves the minimum problem \eqref{OS2}-\eqref{Q_n}
is called the {\it optimal polynomial} for $\mu$ in the class $\mathcal{Q}_n$.
This minimum problem was solved by G. Szeg\H{o}
by showing that the optimal polynomial $p_n(z)$ exists,
is unique and can be expressed in terms of {\it orthogonal polynomials
on the unit circle}
with respect to the measure $\mu$ (see Theorem \ref{Sz2} below).

To state Szeg\H{o}'s solution of the minimum problem \eqref{OS2}-\eqref{Q_n},
we first recall some facts from the theory of orthogonal polynomials on the unit circle (OPUC).

The system of orthogonal polynomials on the unit circle associated with the measure $\mu$:
$$\{\I_n(z) = \I_n(z;f), \q z=e^{i\la}, \q n\in \mathbb{Z}_+\}$$
is uniquely determined by the following two conditions:
\begin{itemize}
\item[(i)] \, $\I_n(z) = \kappa_nz^n + \cdots +l_n$
is a polynomial of degree $n$, in which the leading coefficient
$\kappa_n$ is positive;
\item[(ii)] $(\I_k,\I_j)_\mu = \de_{kj}$ for arbitrary $k,j\in \mathbb{Z}_+$, where
$\de_{kj}$ is the Kronecker delta.
\end{itemize}
Define the monic ($p_n(z)$) and the reciprocal ($p^*_n(z)$) polynomials
(see, e.g., Simon \cite{Sim-1}, p. 2):
\bea
\label{oppp}
&&p_n(z):=p_n(z,\mu)=\kappa^{-1}_n\I_n(z)=z^n + \cdots +l_n\kappa_n^{-1}, \\
\label{opp*}
&&p^*_n(z):= p^*_n(z,\mu)=z^{n}\ol {p_n(1/\ol z)}=\ol l_n\kappa_n^{-1}z^n + \cdots +1.
\eea
We have
\beq
\label{opp}
||p_n||_{2,\mu}=||p^*_n||_{2,\mu}=\kappa_n^{-2}.
\eeq

The polynomials $p_n(z)$ and $p^*_n(z)$ satisfy {\it Szeg\H{o}'s recursion relation}
(see Simon \cite{Sim-1}, p. 56):
\beq
\label{opp1}
p_{n+1}(z)=zp_{n}(z)-\ol v_{n+1} p^*_{n}(z), \q n\in\mathbb{Z}_+
\eeq
where
\beq
\label{opp2}
v_{n+1} = -\ol {p_{n+1}(0)}=\ol l_{n+1}\kappa_{n+1}^{-1},\q |v_{n+1}|<1.
\eeq
In view of \eqref{opp1} we have (see Simon \cite{Sim-1}, p. 56)
\beq
\label{opp3}
||p_n||_{2,\mu}^2=(1-|v_n|^2)||p_{n-1}||_{2,\mu}^2=\prod_{j=1}^n(1-|v_j|^2), \q n\in\mathbb{N}.
\eeq
From \eqref{opp} and \eqref{opp3} we obtain
\beq
\label{opp4}
\kappa_n^{2}\kappa_{n+1}^{-2} =1-|v_{n+1}|^2.
\eeq
The parameters $v_n:=v_n(\mu)$ ($n\in\mathbb{N}$), which play an
important role in the theory of OPUC,  are called {\it Verblunsky's coefficients}
(also known as the Szeg\H{o}, Schur, and canonical moments;
see Simon \cite{Sim-1}, Sect. 1.1, and Dette and Studden \cite{DS}, Sect. 9.4).

\vskip1mm
\n
{\it Note.} The term "Verblunsky coefficient" is from Simon \cite{Sim-1}.
Observe that we write $v_{n+1}$ for Simon's $\al_{n}$, and so one has
$n\in\mathbb{N}$ for Simon's $n\in\mathbb{Z}_+$.
Our notational convention is already established in the time-series
literature and is more convenient in our context of the PACF
(defined below), where $n\in\mathbb{N}$
(see Bingham \cite{Bin1}, Brockwell and Davis \cite{BD}, Sec. 5.2,
Inoue \cite{In2}, Pourahmadi \cite{Po}, Sec. 7.3).

The following result shows that Verblunsky's coefficients provide
a convenient way for the parametrization of probability measures on the unit
circle $\mathbb{T}$ (see, e.g., Verblunsky \cite{Ver1,Ver2},
Ramsey \cite{Ram}, Simon \cite{Sim-1}, p. 2).

\begin{thm}[Verblunsky \cite{Ver1}]
\label{Ver}
Let $\mathbb{D}^\f:=\times_{k=0}^\f\mathbb{D}$ be the set of complex
sequences $v:=(v_n, n\in \mathbb{N})$ with $v_n\in \mathbb{D}$.
The map $\mathcal{S}: \mu \longmapsto v$ is a bijection between the
set of nontrivial probability measures $\{\mu\}$ on $\mathbb{T}$ and $\mathbb{D}^\f$.
\end{thm}

This result was established by Verblunsky \cite{Ver1} in 1935, in connection with OPUC.
It was re-discovered by Ramsey \cite{Ram} in 1974, in connection
with parametrization of time-series models.

\vskip2mm
\n \underline{{\it Partial autocorrelation function(PACF).}}
For a stationary process $X(t)$ with a non-trivial spectral measure $\mu$
the {\it partial autocorrelation function} (PACF) of $X(t)$, denoted by
$\pi_n=\pi_n(\mu)$ ($n\in\mathbb{N}$),
is defined to be the correlation coefficient between the forward and backward
residuals in the linear prediction of the variables $X(n)$ and $X(0)$ on the basis
of the intermediate observations $X(1), \ldots, X(n-1)$, that is,
$$\pi_n:={\rm corr}(X(n)-\widehat X(n), X(0)-\widehat X(0)).$$

It turns out that the Verblunsky coefficients $v_n$ and the PACF $\pi_n$ coincide,
that is, $v_n=\pi_n$ for all  $n\in \mathbb{N}$ (see Dette and Studden \cite{DS}, Sect. 9.6).
Thus, the Verblunsky sequence  $v:=(v_n, n\in \mathbb{N})$ provides a link
between OPUC and time-series analysis, and, in view of the equality $v_n=\pi_n$,
the Verblunsky bijection gives a {\it simple and unconstrained parametrization}
of stationary processes, in contrast to using the covariance function,
which has to be positive-definite.

The next result by Szeg\H{o} solves the minimum problem \eqref{OS2}-\eqref{Q_n}
(see, e.g., Grenander and Szeg\H{o} \cite{GS}, p. 38).

\begin{thm}[Szeg\H{o}]
\label{Sz2}
The unique solution of the minimum problem \eqref{OS2}-\eqref{Q_n}
is the monic polynomial  $p_n(\mu):=p_n(z,\mu)$ 
given by formula \eqref{oppp},
and the minimum in \eqref{OS2} is equal to $||p_n||_{2,\mu}=\kappa^{-2}_n$ (see \eqref{opp}).
\end{thm}

Thus, for the prediction error $\sigma_n^2(\mu)$ we have the following formula:
\bea
\label{mm1}
\si^2_n(\mu)=
\min_{\{q_n\in \mathcal{Q}_n\}}\left\|q_n\right\|^2_{2,\mu}
=\left\|p_n(\mu)\right\|^2_{2,\mu}=\kappa^{-2}_n.
\eea
\begin{rem}
\label{Rem22}
{\rm Define
\beq
\label{Q_n*}
\mathcal{Q}^*_n: =\left\{q_n:  q_n(z)=c_0z^{n} + c_1z^{n-1}+ \cdots c_n, \, c_n=1\right\},
\eeq
and observe that the classes of polynomials $\mathcal{Q}_n$ and $\mathcal{Q}^*_n$
defined in \eqref{Q_n} and \eqref{Q_n*}, respectively, differ by normalization:
in \eqref{Q_n*} we have $ c_n=1$, while in \eqref{Q_n} we have $ c_0=1$.
Also, the optimal polynomial for $\mu$ in the class $\mathcal{Q}^*_n$
is the reciprocal polynomial $p^*_n(z)$ (see \eqref{opp*}).
Taking into account \eqref{opp}, we have the following formula for the prediction
error $\si^2_n(\mu)$  in terms of the optimal polynomial $p^*_n(z)$:
\bea
\label{mmm1}
\si^2_n(\mu)=
\min_{\{q_n\in \mathcal{Q}^*_n\}}\left\|q_n\right\|^2_{2,\mu}
=\left\|p^*_n(\mu)\right\|^2_{2,\mu}.
\eea}
\end{rem}
\begin{rem}
\label{det}
{\rm
Denote by $D_n=D_n(\mu):=\det[r(t-s), \, t,s =0,1,\ldots n]$ the $n^{th}$
Toeplitz determinant generated by the measure $\mu$,
where $r(t)$ is the covariance function given by \eqref{i1}.
Taking into account that $\kappa^{2}_n={D_{n-1}/D_n}$
(see, e.g., Grenander and Szeg\H{o} \cite{GS}, p. 38),
in view of \eqref{mm1} we obtain the following formula
for the prediction error $\sigma_n^2(\mu)$ in terms of $D_n(\mu)$:}
\bea
\label{mm1d}
\si^2_n(\mu)=\frac{D_n(\mu)}{D_{n-1}(\mu)}.
\eea
\end{rem}
\begin{rem}
\label{nn1}
{\rm
In view of \eqref{mm1}, the formulas \eqref{opp3} and \eqref{opp4} can be written as follows
\beq
\label{opp3v}
\si^2_n(\mu)=\prod_{j=1}^n(1-|v_j|^2)\q {\rm and } \q
\frac{\si^2_{n+1}(\mu)}{\si^2_n(\mu)} =1-|v_n|^2.
\eeq
From the second formula in \eqref{opp3v}, it follows that the convergence of the sequences $|v_n|$
and ${\si_{n+1}(\mu)}/{\si_n(\mu)}$ are equivalent, and, the greater
the limiting value of $|v_n|$, the faster the rate of decrease of $\si_n(\mu)$.}
\end{rem}

For a general measure $\mu$ of the form \eqref{di1} the asymptotic relation
\beq
\label{Rv}
\lim_{n\to\f}v_n(\mu)=0
\eeq
is of special interest in the theory of OPUC. In this respect the
following question arises naturally: what is the 'minimal' sufficient condition
on the measure $\mu$ ensuring the relation \eqref{Rv}?
The next result of Rakhmanov \cite{Ra4a,Ra4} shows that for \eqref{Rv}
(or equivalently, for $\lim_{n\to\f}{\si_{n+1}(\mu)}/{\si_n(\mu)}=1$)
it is enough only to have a.e. positiveness on $\mathbb{T}$ of the
spectral density $f$ (see also Babayan et al. \cite{BGT} and Simon \cite{Sim-1}, p. 5).

\begin{thm}[Rakhmanov \cite{Ra4a}]
\label{Rah}
Let the measure $\mu$ have the form \eqref{di1} with $f>0$ a.e. on $\mathbb{T}$.
Then the asymptotic relation \eqref{Rv} is satisfied.
\end{thm}
Note that the converse of Rakhmanov's theorem, in general,
is not true. A partial converse of Rakhmanov's theorem
is stated in Theorem \ref{cG}.

Bello and L\'opez \cite{BL} proved the following extension of Rakhmanov's theorem:
{\it Let $\G_\de$ be a closed arc of the unit circle of length
$2\de$ $(0<\de\leq\pi)$, and let $\mu$ and $(v_n)$ be as in Rakhmanov's Theorem.
Assume that the measure $\mu$ is supported on the arc $\G_\de$ with $f>0$
a.e. on $\G_\de$. Then $\lim_{n\to\f}|v_n|=\cos(\de/2)$.}
The case $\de=\pi$ corresponds to Rakhmanov's theorem.

\sn{Properties of the prediction error $\si_{n}^2(f)$}

In what follows we assume that the spectral measure $\mu$ is absolutely
continuous with spectral density $f$, and instead of $\si^2_n(\mu)$, $p_n(\mu)$ and $D_n(\mu)$
we use the notation $\si^2_n(f)$, $p_n(f)$ and $D_n(f)$, respectively.

In the next proposition we list a number of properties of the prediction error $\sigma_n^2(f)$.
The proof can be found in Babayan and Ginovyan \cite{BG}.
\begin{pp}
\label{pp3}
The prediction error $\sigma_n^2(f)$ possesses the following properties.
\begin{itemize}
\item[(a)]
$\sigma_n^2(f)$ is a non-decreasing functional of $f$:
$\sigma_{n}^2(f_1)\leq\sigma_n^2(f_2)$ when $f_1(\lambda)\leq f_2(\lambda)$, $\la\in [-\pi,\pi]$.
\item[(b)]
If $f(\la)=g(\la)$ almost everywhere on $[-\pi,\pi]$, then $\si_n^2(f) = \si_n^2(g).$
\item[(c)]
For any positive constant $c$ we have $\si_n^2(cf) = c\si_n^2(f)$.
\item[(d)]
If $\bar f(\la)=f(\la-\la_0)$, $\la_0\in [-\pi,\pi]$, then $\si_n^2(\bar f) = \si_n^2(f).$
\end{itemize}
\end{pp}

\s {Asymptotic behavior of the prediction error $\de_n(f)$
for nondeterministic processes}
\label{ndp}
In this section we study the asymptotic behavior of the the relative
prediction error $\de_n(f) = \si^2_n(f) - \si^2(f)$
for nondeterministic processes, and review some important known results.

\sn{Asymptotic behavior of $\de_n(f)$ for short-memory processes}

\n Recall that a {\sl short memory\/} processes is a second order stationary
processes possessing a spectral density $f$ which is bounded away from
zero and infinity.

\ssn{Exponential rate of decrease of $\de_n(f)$} 

\n The first result of this type goes back to the Grenander and Rosenblatt
\cite{GR1}. The next theorem was proved by Ibragimov \cite{I-2}.
\begin{thm}[Ibragimov \cite{I-2}]
\label{I1}
A necessary and sufficient condition for
\beq
\label{eex1}
\de_n(f) = O(q^n), \q q=e^{- c}, \q c > 0, \q n \ra \f
\eeq
is that $f(\la)$ is a spectral density of a short-memory
process, and $1/{f(\la)} \in A_c$, where $A_c$ is the class of
$2\pi$--periodic continuous functions $\I(\la), \, \la\in\mathbb{R}$,
admitting an analytic continuation into the strip
$z=\la+i\mu$, $-\f<\la<\f$, $|\mu|\le c$.
\end{thm}

\n Observe that (\ref{eex1}) will be true for all $c>0$ if and only if the
analytic continuation of $f(\la)$ is an entire function of $z=\la+i\mu$.

\n Thus, to have exponential rate of decrease to zero for $\de_n(f)$ the
underlying model should be short-memory process with sufficiently smooth
spectral density.

\ssn {Hyperbolic rate of decrease of $\de_n(f)$}
Here we are interested in estimates for $\de_n(f)$ of type
\bea
\label{eex2}
\de_n(f) &=& O(n^{- \g}), \q \g > 0, \q n \ra \f.\\
\label{eex3}
\de_n(f) &=& o(n^{- \g}), \q \g > 0, \q n \ra \f.
\eea

Bounds of type (\ref{eex2}) with $\g > 1$ for different classes of
spectral densities were obtained by Baxter \cite{Bax}, Devinatz \cite{Dev1},
Geronimus \cite{Ger-1, Ger-2}, Grenander and Rosenblatt \cite{GR1}, Grenander and
Szeg\H{o} \cite{GS}, and others (see, e.g., Ginovyan \cite{G-4}, and references
therein).
\n The most general result in this direction has been obtained by
Ibragimov \cite{I-2}. To state Ibragimov's theorem, we first introduce
H\"older classes of functions.

For a function $\I (\la)\in L^p(\mathbb{T})$, we define its
$L^p$-modulus of continuity by
\bea
\label{eex4}
\om_p(\I; \de) = \sup_{0<|t|\le\de}||\I (\cd + t) - \I (\cd)||_p,
\q \de > 0.
\eea
Given numbers $0 < \al < 1$, $r \in \mathbb{Z}_0:= \{0,1,2\ldots\}$,
and $p \ge 1$, we put $\g: = r + \al$.
A H\"older class of functions, denoted by $H_p(\g)$, is defined
to be the set of those functions $\I(\la) \in L^p(\mathbb{T})$ that
have $r$-th derivative $\I^{(r)}(\la)$, such that
$\I^{(r)}(\la) \in L^p(\mathbb{T})$ and $\om_p(\I^{(r)}; \de) = O(\de^\al)$ as $\de\to0$.

\begin{thm}[Ibragimov \cite{I-2}]
\label{I2}
A necessary and sufficient condition for
\bea
\label{eex6}
\de_n(f) =O(n^{- \g}), \q \g = 2(r + \al) > 1; \,\,
0 < \al < 1, \, r \in \mathbb{Z}_0, \q {\rm as} \q n \ra \f
\eea
\n is that $f(\la)$ is a spectral density of a short-memory
process belonging to $H_2(\g)$.
\end{thm}

Bounds of type (\ref{eex3}) for short memory models have been obtained
by Baxter \cite{Bax}, Devinatz \cite{Dev1}, Hirshman \cite{Hir} and
Golinskii \cite{Gol-1}.
\n For 'large' values of $\g $ ($\g > 1$), Hirshman has obtained the
following necessary and sufficient condition for (\ref{eex3})
(see Hirshman \cite{Hir}, p. 314).

\begin{thm}[Hirshman \cite{Hir}]
\label{Hir}
If $X(t)$ is a short-memory process, then
$\de_n(f) =o(n^{-\g})$ with $\g > 1$ as $n \ra \f$
if and only if
$\sum_{|t|\ge n}{|r(t)|^2} =o(n^{-\g})$ as $n \ra \f$,
where $r(t)$ is the covariance function of $X(t)$.
\end{thm}

The next theorem was proved by G. Baxter (see Baxter \cite{Bax},
Theorem 3.1).
\begin{thm}[Baxter \cite{Bax}]
\label{Bax1}
If $X(t)$ is a short-memory process with covariance function $r(t)$
satisfying $\sum_{t=1}^\f{t^\g|r(t)|} <\f$, $\g > 0$, then
$\de_n(f) =o(n^{-2\g})$ as $n \ra \f$.
\end{thm}

\begin{rem}
{\rm It follows from Theorem \ref{I2} that if $\de_n(f) =O(n^{- \g})$
with $\g > 1$, then the underlying model $X(t)$ is necessarily a short-memory
process. Moreover, as it was pointed out by Grenander and Rosenblatt \cite{GR1}
(see, also, Devinatz \cite{Dev1}, p. 118), if the model is not a short-memory
process, that is, the spectral density $f$ has zeros or is unbounded,
then, in general, we cannot expect $\de_n(f)$  to go to zero faster than $1/n$
as $n \to \f$. This question we discuss in the next section.}
\end{rem}

\sn{Asymptotic behavior of $\de_n(f)$ for anti-persistent and
long-memory processes}

\n In this section we describe the rate of decrease of the relative prediction
error $\delta_n(f)$ to zero as $n \to \infty,$ in the case where
the underlying process $X(t)$ is nondeterministic and is anti-persistent or
has long-memory, that is, the spectral density $f(\la)$ of $X(t)$ has zeros
or is unbounded at a finite number of points, such that
$\ln f(\la) \in L^1(\mathbb{T})$.
This case is of great interest because in many applications the model
is described by such type processes.

\ssn{An example}
We start with an example, which shows that the asymptotic behavior
of the prediction error $\de_n(f)$ essentially depends on the
dependence (memory) structure of the underlying model $X(t)$
(see, e.g., Grenander and Szeg\H{o} \cite{GS}, p. 191).

\begin{exa}
{\rm Consider a first-order moving-average $MA(1)$ process $X(t)$:
$$
X(t)=\vs (t)-b \cd \vs (t-1),\q \vs (t)\sim WN(0,\si^2_\vs),
\q 0\le b\le1.
$$
\n The spectral density is
$f(\la)=\frac{\si^2_\vs}{2\pi}\cd |1-be^{i\la}|^2$
(see formula \eqref{arma}).

a) First assume that $X(t)$ has short-memory, that is, $0\le b<1$.
It is easy to check that 
$$\de_n(f) = \si_n^2(f)-\si^2(f)=\frac{b^{2n}(b^2-b^{4})}{1-b^{2n+2}}
\sim b^{2n}\q {\rm as}\q n\to\f,$$
showing that in this case $\de_n(f)$ goes to zero with exponential rate.

b) Now let $b=1$. We have $f(\la)= \frac{\si^2_\vs}{2\pi}\cd|1-e^{i\la}|^2,$
that is, the process $X(t)$ is anti-persistent. In this case we have
$$\de_n(f) = \si_n^2(f)-\si^2(f)=\frac1{n+2}\sim \frac1{n}\q {\rm as}\q n\to\f,$$
showing that $\de_n(f)$  goes to zero fat precisely the rate $1/n$.
The slow rate is due to the presence of a zero of $f(\la)$ at $\la=0$
(see Grenander and Szeg\H{o} \cite{GS}, p. 191).

It can be shown that for models with spectral densities of the form
$
f(\la) = \frac{\si^2_\vs}{2\pi}\cd|q(e^{i\la})|^2,
$
where $q(e^{i\la})$ is a polynomial with at least one root on the unit
circle $\mathbb{T}$, we have $\de_n(f) \sim\frac{1}n$ as $n \to \f.$
}
\end{exa}

\ssn {The ARFIMA$(p,d,q)$ Model}

As it was mentioned in Section \ref{mem} a well-known example of
processes that displays long memory or is anti-persistent is an
ARFIMA$(p,d,q)$ process $X(t)$  with spectral density
$f_X(\la)$ given by (see \eqref{AA}):
\begin{equation}
\label{lm3}
f_X(\la)=|1-e^{-i\la}|^{-2d}f(\la), \q d<1/2,
\end{equation}
where $f(\la)$ is the spectral density of an ARMA$(p,q)$ process,
given by (\ref{arma}).
Recall that for $ 0<d<1/2$ the model $X(t)$ specified by spectral density (\ref{lm3})
displays long-memory, for $d<0$ it is anti-persistent,
and for $d=0$ it displays short-memory.

The following theorem was proved by A. Inoue (see Inoue \cite{In2}, Theorem 4.3).
\begin{thm}[Inoue \cite{In2}]
\label{TLM1}
Let $f_X(\la)$ have the form (\ref{lm3}) with $ 0<d<1/2$, where $f(\la)$
is the spectral density of an ARMA$(p,q)$ process. Then
\beq
\label{lm4}
\de_n(f) \sim\frac{d^2}n \q {\rm as} \q n \to \f.
\eeq
\end{thm}
\begin{rem}
\label{lg1}
{\rm Note that for ARFIMA$(0,d,0)$ model the asymptotic relation (\ref{lm4})
remains valid for all $d<1/2$ ($d\neq 0$). In this case, for the Verblunsky
coefficients (parameters)
$v_n$ we have $v_n =\frac{d}{n-d+1}$ (see Golinskii \cite{Gol-1}, p. 703).}
\end{rem}

\ssn {The Jacobian Model}

Another well-known example of processes that displays long memory
or is anti-persistent is the Jacobian model.
We say that a stationary process $X(t)$ is a Jacobian process,
and the corresponding model is a Jacobian model, if its spectral density
$f(\la)$ has the following form:
\beq
\label{J1}
f(\la) = f_1(\la)\prod_{k=1}^m|e^{i\la}-
e^{i\la_k}|^{-2d_k},
\eeq
where $f_1(\la)$ is the spectral density of a short-memory process,
the points $\la_k\in[-\pi, \pi]$ are distinct, and $d_k\leq1/2$, $k=1,\ldots,m$.

The asymptotic behavior of $\de_n(f)$ as $n\to\f$ for Jacobian model \eqref{J1}
has been considered in a number of papers (see, e.g., Golinskii \cite{Gol-1},
Grenander and Rosenblatt \cite{GR1}, Ibragimov \cite{I-2}, Ibragimov and
Solev \cite{IS}.)

The following theorem was proved in Ibragimov and Solev \cite{IS}.
\begin{thm}[Ibragimov and Solev \cite{IS}]
\label{IS}
Let $f(\la)$ have the form \eqref{J1}, where $f_1(\la)$ is the spectral density
of a short-memory process, the points $\la_k\in[-\pi, \pi]$ are distinct,
and $d_k\leq 1/2$, $k=1,\ldots,m$. If the function $f_1(\la)$ satisfies a
Lipschitz condition with exponent $\al\ge1/2$, then
\beq
\label{IbS}
\de_n(f) \sim{1}/n \q {\rm as} \q n \to \f.
\eeq
\end{thm}

\section{Asymptotic behavior of the finite prediction error $\si_n^2(f)$
for deterministic (singular) processes} \label{d}

\sn {Background}

The linear prediction problem has been studied most intensively
for nondeterministic processes, that is, in the case where the
prediction error is known to be positive $(\si^2(f) >0)$
(see Section \ref{ndp}).

In this section we focus on the less investigated case - deterministic
processes, that is, when  $\si^2(f) =0$.
This case is not only of theoretical interest, but is also important
from the point of view of applications.
For example, as pointed out by Rosenblatt \cite{Ros} (see also Pierson \cite{Pl}),
situations of this type arise in Neumann's theoretical model of
storm-generated ocean waves. Such models are also of interest
in meteorology (see, e.g., Fortus \cite{For}).

Only few works are devoted to the study of the speed of convergence of
$\de_n(f)=\si^2_n(f)$ to zero as $n\to\infty$, that is, the asymptotic
behavior of the prediction error for deterministic processes.
One needs to go back to the classical work of M. Rosenblatt \cite{Ros}.
Using the technique of orthogonal polynomials on the unit circle,
M. Rosenblatt investigated the asymptotic behavior of the prediction
error $\si^2_n(f)$ for deterministic processes in the following two cases:
\begin{itemize}
\item[(a)]
the spectral density $f(\la)$ is continuous and positive on a
segment of $[-\pi,\pi]$ and zero elsewhere.

\item[(b)]
the spectral density $f(\la)$ has a very high order of contact with zero
at points $\la=0, \pm\pi$, and is strictly positive otherwise.
\end{itemize}

Later the problems (a) and (b) were studied by Babayan  \cite{Bb-1,Bb-2},
Babayan and Ginovyan \cite{BG-1,BG,BG-2},
and Babayan et al. \cite{BGT}
(see also Davisson \cite{Dav1} and Fortus \cite{For}), where some generalizations and
extensions of Rosenblatt's results have been obtained.

\medskip
We start by describing Rosenblatt's results concerning the asymptotic behavior
of the prediction error $\si^2_n(f)$, obtained in Rosenblatt \cite{Ros} for the
above stated cases (a) and (b).

\sn {Rosenblatt's results about speed of convergence}

For the case (a) above, M. Rosenblatt proved in \cite{Ros}
that the prediction error $\si^2_n(f)$ decreases to zero exponentially
as $n\to\f$. More precisely, M. Rosenblatt proved
the following theorem.

\begin{thm}[Rosenblatt's first theorem]
\label{R1}
Let the spectral density $f$ of a discrete-time stationary
process $X(t)$ be positive and continuous on the
segment $[\pi/2-\alpha, \pi/2+\alpha],$ $0<\alpha<\pi,$
and zero elsewhere. Then the
prediction error $\si^2_n(f)$ approaches zero
exponentially as $n\to\f$.
More precisely, the following asymptotic relation holds:
\beq
\label{nd1}
\si^2_n(f)\simeq\left(\sin(\alpha/2)\right)^{2n+1}
\q {\rm as} \q n\to\f.
\eeq
\end{thm}
Thus, when the spectral density $f$ is continuous and vanishes on an
entire segment, then the prediction error $\si^2_n(f)$ approaches zero
with a sufficiently high speed, namely as a geometric progression with
common ratio $\sin^2(\alpha/2)<1$.
Notice that \eqref{nd1} implies that
\beq
\label{nd2}
\lim_{n\to\f}\sqrt[n]{\sigma_n(f)} =\sin(\alpha/2).
\eeq

Concerning the case (b) above, for a specific deterministic process $X(t)$, Rosenblatt proved in \cite{Ros} that
the prediction error $\si^2_n(f)$ decreases to zero {\it like a power} as $n\to\f$.
More precisely, the deterministic process $X(t)$ considered in
Rosenblatt \cite{Ros} has the spectral density
\beq
\label{nd4}
f_a(\la):=\frac{e^{(2\la-\pi)\varphi(\la)}}{\cosh\left(\pi\varphi(\la)\right)},
\q f_a(-\la)=f_a(\la),\q 0\leq\la\leq\pi,
\eeq
where $\varphi(\la)=(a/2)\cot\la$ and $a$ is a positive parameter.

Using the technique of orthogonal polynomials on the unit circle
and Szeg\H{o}'s results, Rosenblatt \cite{Ros} proved the following theorem.
\begin{thm}[Rosenblatt's second theorem]
\label{R2}
Suppose that the process $X(t)$ has spectral density $f_a$ given by \eqref{nd4}.
Then the following asymptotic relation for the prediction error $\si^2_n(f_a)$
holds:
\beq
\label{nd6}
\si^2_n(f_a)\sim\frac{\Gamma^2\left(({a+1)}/2\right)}
{\pi 2^{2-a}} \ n^{-a}
\q {\rm as} \q n\to\f.
\eeq
\end{thm}
Note that the function in \eqref{nd4} was first considered by Pollaczek \cite{Pol}, and then by Szeg\H{o} \cite{S1}, as a weight-function of a class of orthogonal polynomials that serve as illustrations for certain 'irregular' phenomena in the theory of orthogonal polynomials.
For the function $f_a$ in \eqref{nd4}, we have the following asymptotic
relation (for details see Szeg\H{o} \cite{S1},
{Babayan and Ginovyan \cite{BG}}, and Section \ref{dex3}):
\beq \label{tt2}
f_a(\la)\sim
\left \{
\begin{array}{ll}
 2e^a\exp\left\{-{a\pi}/{|\la|}\right\} & \mbox{as $\la\to0$},\\
2\exp\left\{-{a\pi}/{(\pi-|\la|)}\right\} & \mbox{as $\la\to\pm\pi$}.
\end{array}
\right.
\eeq
Thus, the function $f_a$ in \eqref{nd4} has a very high order of contact
with zero at points $\la=0,\pm\pi$, due to which the process with spectral
density $f_a$ is deterministic and the prediction error $\si^2_n(f_a)$
in \eqref{nd6} decreases to zero like a power as $n\to\f$.

\begin{rem}
{\rm
{In view of formulas in \eqref{opp3v}, under the conditions of
Theorem \ref{R1}, we have $$\lim_{n\to\f}{\si^2_{n+1}(f)}/{\si^2_n(f)}=\sin^2(\al/2) \q {\rm and } \q
\lim_{n\to\f}|v_n(f)|=\cos(\al/2).$$
Similarly, under the conditions of Theorem \ref{R2},
we have
$$\lim_{n\to\f}{\si_{n+1}(f_a)}/{\si_n(f_a)}=1 \q {\rm and } \q
\lim_{n\to\f}v_n(f_a)=0,$$
where $v_n(f)$ and $v_n(f_a)$ are the
Verblunsky coefficients corresponding to functions $f$ and $f_a$,
respectively.}}
\end{rem}

In the papers Babayan \cite{Bb-1,Bb-2}, Babayan and Ginovyan \cite{BG-1,BG,BG-2},
and Babayan et al. \cite{BGT}, the above stated Theorems \ref{R1} and \ref{R2}
were extended to broader classes of spectral densities.

Concerning Theorem \ref{R1}, was described an extension of the
asymptotic relation \eqref{nd2}
to the case of several arcs, without having to stipulate continuity
of the spectral density $f$.

As for the extension of Theorem \ref{R2}, it was proved that if the spectral
density $f$ is such that the sequence $\si_n(f)$ is weakly varying
(a term defined in Section \ref{ww}) and if, in addition, $g$ is
a nonnegative function that can have arbitrary power type singularities, then the sequences $\{\si_n(fg)\}$ and $\{\si_n(f)\}$
have the same asymptotic behavior as $n\to\f$, up to some positive multiplicative factor.
This allows to derive the asymptotic behavior of $\si_n(fg)$ from that of $\si_n(f)$.

Using this result, Rosenblatt's Theorem \ref{R2} was extended in
Babayan and Ginovyan \cite{BG} and in Babayan et al. \cite{BGT} to a
class of spectral densities of the form $f=f_ag$,
where $f_a$ is as in (\ref{nd4}) and $g$ is a nonnegative function that
can have arbitrary power type singularities
(see Corollary \ref{c4.2} in Section \ref{fdR}).

\sn{Extensions of Rosenblatt's first theorem}
\label{d1}

In this section we extend Rosenblatt's first theorem (Theorem \ref{R1})
to a broader class of deterministic processes, possessing spectral densities.
More precisely, we extend the asymptotic relation
\eqref{nd2} to the case of several arcs, without having to stipulate continuity
of the spectral density $f$.
Besides, we state necessary as well as sufficient conditions for the exponential
decay of the prediction error $\si_{n}(f)$.
Also, we calculate the transfinite diameters of some subsets of the unit circle,
and thus, obtain explicit asymptotic relations for
$\si_n(f)$ similar to Rosenblatt's relation \eqref{nd2}.

\ssn{Extensions of Theorem \ref{R1}}

In what follows, by $E_f$ we denote the spectrum of the process $X(t)$, that is,
\beq
\label{Sp}
E_f:=\{e^{i\la}: f(\la)> 0\}.
\eeq
Thus, the closure $\ol E_f$ of $E_f$ is the support of the spectral density $f$.

The next theorem extends Rosenblatt's first theorem (Theorem \ref{R1}).
More precisely, the result that follows extends the asymptotic
relation \eqref{nd2} to the case of several arcs, without having
to stipulate continuity of the spectral density $f$.

\begin{thm}[Babayan et al. \cite{BGT}]
\label{BB1}
Let the support $\ol E_f$ of the spectral density $f$ of the process $X(t)$
consists of a finite number of closed arcs of the unit circle $\mathbb{T}$,
and let $f>0$ a.e. on $\ol E_f$.
Then the sequence $\sqrt[n]{\sigma_n(f)}$ converges, and
\beq
\label{bb2}
\lim_{n\to\f}\sqrt[n]{\sigma_n(f)}  =\tau_f, 
\eeq
where $\tau_f:=\tau(\ol E_f)$ is the transfinite diameter of $\ol E_f$.
\end{thm}
\begin{rem}
\label{rm3.5}
{\rm
In Theorem \ref{R1},
$\ol E_f=\{e^{i\la}: \la\in [\pi/2-\alpha, \pi/2+\alpha]\},$
which represents a closed arc of length $2\al$, and, according to Proposition \ref{pp1}(d),
we have $\tau(\ol E_f)=\sin(2\alpha/4)=\sin(\alpha/2)$.
Thus, the asymptotic relation \eqref{nd2} is a special case of \eqref{bb2}.}
\end{rem}

We will need the following definition, which characterizes
the rate of variation of a sequence of non-negative numbers compared with a geometric progression
(see also Simon \cite{Sim-1}, p. 91).
\begin{den}
\label{ekd1}
(a) A sequence $\{a_n\geq 0, \, n \in\mathbb{N}\}$
is called exponentially neutral if
$$
\lim_{n \rightarrow \infty} \sqrt[n]{a_n} =1.
$$
(b) A sequence $\{b_n\geq 0, \, n \in\mathbb{N}\}$
is called exponentially decreasing if
$$
\limsup_{n \rightarrow \infty} \sqrt[n]{b_n} <1.
$$
\end{den}
For instance, the sequence $\{a_n=n^\al, \, \al\in\mathbb{R}, \,
n \in\mathbb{N}\}$ is exponentially neutral because
$\ln \sqrt[n]{n^\alpha} = ({\alpha}/{n}) \ln {n}\rightarrow 0$ as $n\to\f$.
The geometric progression $\{b_n=q^n, \, 0<q<1, \, n \in\mathbb{N}\}$
is exponentially decreasing because
$\sqrt[n]{b_n} = q^{n/n} = q<1$.
The sequence $\{b_n=n^\al q^n, \, \al\in\mathbb{R}, \, 0<q<1, \,
n \in\mathbb{N}\}$ is also exponentially decreasing because
$\sqrt[n]{b_n} = n^{\alpha/n}q \rightarrow q<1.$
In fact, it can easily be shown that a sequence $\{c_n\geq0, \, n \in\mathbb{N}\}$
is exponentially decreasing if and only if there exists a number
$q$ ($0<q<1$) such that $c_n =O(q^{n})$ as $n\to\f$.

\begin{rem}
\label{RMT1}
{\rm
It follows from relation \eqref{bb2} that the question of exponential
decay of the prediction error $\sigma_n(f)$ as $n\to\f$ is determined
solely by the value of the transfinite diameter of the support $\ol E_f$
of the spectral density $f$, and does not depend on the values of $f$ on $\ol E_f$.
Denoting $\g_n:=\si_n(f)/\tau_f^n $, from \eqref{bb2} we infer that
$\lim_{n \to \infty} \sqrt[n]{\g_n} =1$ and
\begin{equation}
\label{sfg}
\si_n(f)=\tau_f^n\cd \g_n.
\end{equation}
Thus, in the case where $\tau_f<1$, the prediction error $\sigma_n(f)$
is decomposed into a product of two factors, one of which ($\tau_f^n$)
is a geometric progression, and the second ($\g_n$) is an exponentially
neutral sequence.
Also, if $g$ is another spectral density satisfying the conditions of
Theorem \ref{BB1}, then in view of \eqref{sfg}, we have
\begin{equation}
\nonumber
\frac{\si_n(g)}{\si_n(f)}=\left(\frac{\tau_g}{\tau_f}\right)^n\cd\g'_n,
\end{equation}
where $\g'_n$ is an exponentially neutral sequence. It is worth noting
that the last relation does not depend on the structures of the
supports $\ol E_f$ and $\ol E_g$ (viz., the number and the lengths
of arcs constituting these sets, as well as, their location on the unit
circle $\mathbb{T}$).}
\end{rem}

The following result provides a sufficient condition for the
exponential decay of the prediction error $\si_{n}(f)$.
\begin{thm}[Babayan et al. \cite{BGT}]
\label{SC}
If the spectral density $f$ of the process $X(t)$ vanishes on an arc,
then the prediction error $\si_{n}(f)$ decreases to zero exponentially.
More precisely, if $f$ vanishes on an arc $\G_\de\subset \mathbb{T}$
of length $2\de$ $(0<\de<\pi)$, then
\begin{equation}
\label{m2}
\limsup_{n \rightarrow \infty} \sqrt[n]{\sigma_n(f)} \leq \cos({\de}/{2})<1.
\end{equation}
\end{thm}

The next result gives a necessary condition for the exponential decay of $\si_{n}(f)$.
\begin{thm}[Babayan et al. \cite{BGT}]
\label{NC}
A necessary condition for the prediction error $\si_{n}(f)$ to tend to zero exponentially
is that the spectral density $f$ should vanish on a set of positive Lebesgue measure.
\end{thm}

\begin{rem}
\label{RMT3}
{\rm
Theorem \ref{NC} shows that if the spectral density $f$
is almost everywhere positive, then it is impossible to obtain
exponential decay of the prediction error $\si_{n}(f)$,
no matter how high the orders of the zeros of $f$.}
\end{rem}
In view of \eqref{opp3v}, as a consequence of Theorem \ref{BB1},
we obtain the following result.
\begin{thm}[Babayan et al. \cite{BGT}]
\label{CBB1}
Let the support $\ol E_f$ and the spectral density $f$ satisfy the
conditions of Theorem \ref{BB1}. If the sequence of Verblunsky
coefficients $v_n(f)$ converges in modulus, then
\beq
\label{vvv1}
\lim_{n\to\f}|v_n(f)|  =\sqrt{1-\tau_f^2}.
\eeq
\end{thm}
\begin{rem}
\label{vvv}
{\rm It is well-known that for an arbitrary sequence of positive numbers $a_n$
the convergence $a_{n+1}/a_n\to a$ implies the convergence
$\sqrt[n]{a_n}\to a$. The converse, in general, is not true,
that is, the sequence $\sqrt[n]{a_n}$ can be convergent,
while $a_{n+1}/a_n$ divergent.  Indeed, for the sequence $a_n$:
\begin{equation}
\nonumber
a_n: = \left \{
\begin{array}{ll}
2^{-3k}& \mbox{if \, $n=2k-1$}\\
2^{-(3k+1)}& \mbox{if \, $n=2k$},
           \end{array} \quad k\in\mathbb{N},
           \right.
\end{equation}
we have $\lim_{n\to\f}\sqrt[n]{a_n}=2^{-3/2}$,
while the limit $\lim_{n\to\f}a_{n+1}/a_n$ does not exists.

In the context of the considered sequences, $|v_n(f)|$ and $\sqrt[n]{\sigma_n(f)}$,
in view of \eqref{opp3v}, we can assert that the convergence of $|v_n(f)|$ (or equivalently the convergence of ${\si_{n+1}(f)}/{\si_n(f)}$)
implies the convergence of $\sqrt[n]{\sigma_n(f)}$, but not the converse.
Hence, the condition of convergence (in modulus) of Verblunsky's sequence 
in Theorem \ref{CBB1} is essential.}
\end{rem}

As a consequence of Theorem \ref {BB1} we obtain the following result
(see Geronimus \cite{Ger}), which is a partial converse of Rakhmanov's theorem:
\begin{thm}
\label{cG}
If the sequence $\sigma_n(f)$ satisfies the following condition:
\beq
\label{gt}
\limsup_{n \rightarrow \infty} \sqrt[n]{\sigma_n(f)}=1
\eeq
(in particular, if $\lim_{n\to\f}v_n(f)=0$),
then $\ol E_f=\mathbb{T}$, i.e.
the spectrum of the process is dense in $\mathbb{T}$.
\end{thm}

The next theorem extends Theorem \ref{BB1}.
\begin{thm}[Babayan \cite{Bb-1}]
\label{B3}
Let $E_f$ be the spectrum of a stationary process $X(t)$ possessing
a spectral density $f(\la)$, that is, $E_f=\{\la: f(\la)>0\}$,
and let $\tau(E_f)$, $\tau_*(E_f)$ and $\tau^*(E_f)$ be the transfinite
diameter and the inner and the outer transfinite diameters of $E_f$,
respectively (for definition of $\tau_*(E_f)$ and $\tau^*(E_f)$
see formula \eqref{tdo1}).
Then the following assertions hold.
\begin{itemize}
\item[(a)]
The following inequalities hold:
\begin{align}
\label{tt21}
&\limsup_{n\to\f}(\si_{n}(f))^{1/n} \le\tau^*(E_f),\\
\label{tt25}
&\liminf_{n\to\f}(\si_{n}(f))^{1/n} \ge \tau_*(E_f).
\end{align}
\item[(b)]
If the spectrum $E_f$ consists of a countable number of open arcs of the unit circle $\mathbb{T}$ and is $\tau$-measurable, that is, $\tau_*(E_f) =\tau^*(E_f)=\tau(E_f)$, then
\beq
\label{23g}
\lim_{n \rightarrow \infty} \sqrt[n]{\sigma_n(f)}=\tau(E_f).
\eeq
\end{itemize}
\end{thm}
As an immediate corollary of Theorem \ref{B3} we have the following result.
\begin{cor}
\label{br7}
A sufficient condition for the
prediction error $\si_{n}(f)$ of a deterministic stationary sequence to
decrease to zero at least exponentially as $n\to\f$, that is,
$\sigma_n(f) =O(e^{-b n})$ for some $b>0$, is that the outer transfinite
diameter of
the spectrum $E_f$ should be less than unity.
\end{cor}

\ssn{Examples. Calculation of transfinite diameters of some
special sets}
\label{EX}
Motivated by Theorems \ref{R1} and \ref{BB1} and Remark \ref{rm3.5},
the following question
arises naturally: calculate the transfinite diameter $\tau(\ol E_f)$ 
of the set $\ol E_f$ consisting of several closed arcs of the unit
circle $\mathbb{T}$, and thus, obtain an explicit
asymptotic relation for the prediction error $\si_n(f)$ similar to Rosenblatt's
relation \eqref{nd2}.
The calculation of the transfinite diameter
(and hence, the capacity and the Chebyshev constant) is a challenging problem
(for details see Section \ref{DCC}), and in
only very few cases has the transfinite diameter been exactly calculated (see, e.g., Landkof \cite{La}, p. 172-173, Ransford \cite{Ran}, p.135,
and Proposition \ref{pp1}).
One such example provides Theorem \ref{R1}, in which case the transfinite diameter of the set
$\ol E_f:=\{e^{i\la}: \la\in [\pi/2-\alpha, \pi/2+\alpha]\}$
is $\sin(\alpha/2)$.
Below we give some other examples, where we can explicitly calculate
the Chebyshev constant (and hence the transfinite diameter and the
capacity) by using some properties of the transfinite diameter, stated
in Proposition \ref{pp1}, and results due to Fekete \cite{F}
and Robinson \cite{Rob} concerning the relation between the
transfinite diameters of related sets (see Propositions \ref{FT}
and \ref{RT2}).

The examples given below show that Fekete's formula \eqref{F0} and Robinson's formula \eqref{F10} give a
simple way to calculate the transfinite diameters of some
subsets of the unit circle, based only on the formula
of the transfinite diameter of a line segment (see Proposition \ref{pp1}(e)).

\begin{figure}[ht]%
\centering
\includegraphics[width=0.8\textwidth]{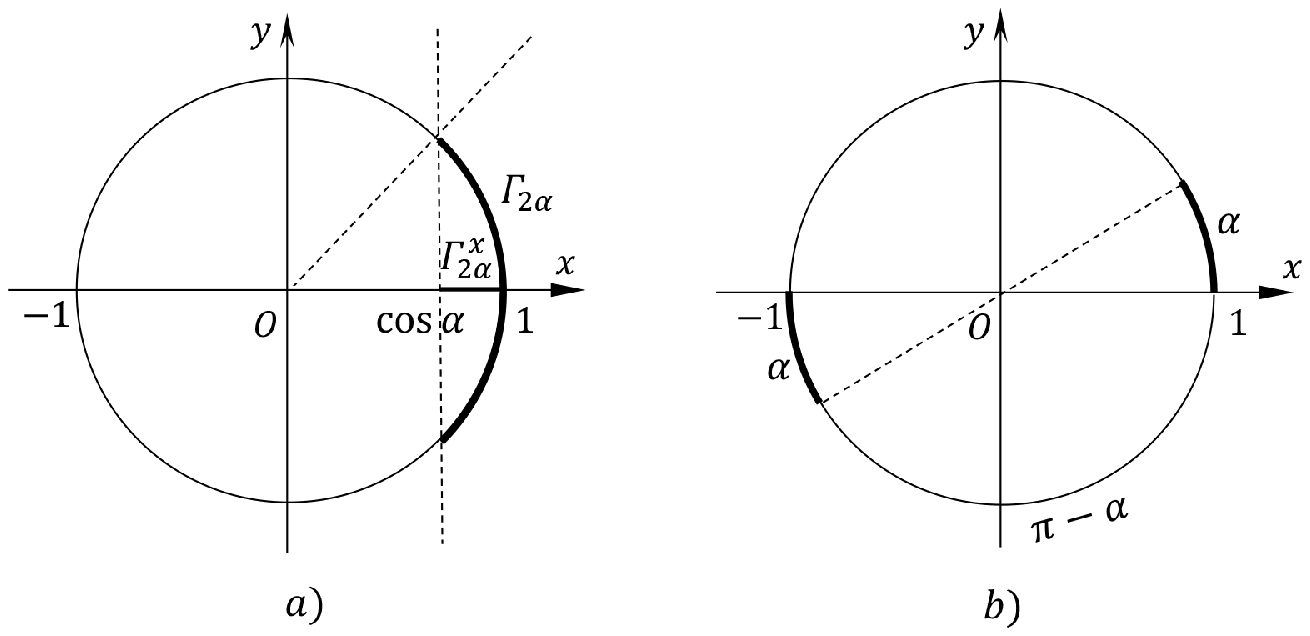}
\caption{a) The sets $\G_{2\al}$ and $\G^x_{2\al}$.
b) The set $\G(k,\al)$ with $k=2$.}
\label{fig1}
\end{figure}

We will use the following notation:
given $0<\be<2\pi$ and $z_0=e^{i\theta_0}$, $\theta_0\in[-\pi,\pi)$,
we denote by $\G_\be(\theta_0)$
an arc of the unit circle of length $\be$ which is symmetric
with respect to the point $z_0=e^{i\theta_0}$: 
\beq
\label{arc}
\G_\be(\theta_0):=\{e^{i\theta}: \, |\theta-\theta_0|\leq\be/2\}
=\{e^{i\theta}: \, \theta\in[\theta_0-\be/2, \theta_0+\be/2]\}.
\eeq
\begin{exa}
\label{ex1}
{\rm Let $\G_{2\al}:=\G_{2\al}(0)$. Then the projection $\G_{2\al}^x$
of $\G_{2\al}$ onto the real axis is the segment $[\cos\al,1]$
(see Figure \ref{fig1}a)),
and by Proposition \ref{pp1}(e) for the transfinite diameter $\tau(\G_{2\al}^x)$ we have
$$\tau(\G_{2\al}^x) =\frac{1-\cos\al}4=\frac{\sin^2(\al/2)}2.$$
Hence, according to Robinson's formula \eqref{F10}, we obtain
\beq
\label{F17}
\tau(\G_{2\al}) =[2\tau(\G_{2\al}^x)]^{1/2}=\left[2\frac{\sin^2(\al/2)}2\right]^{1/2} =\sin(\al/2).
\eeq
Taking into account that the transfinite diameter is invariant with respect
to rotation (see Proposition \ref{pp1}(b)), from \eqref{F17} for any
$\theta_0\in[-\pi,\pi)$ we have
\beq
\label{F17a}
\tau(\G_{2\al}(\theta_0)) = \sin(\al/2).
\eeq}
\end{exa}

\begin{rem}
{\rm Notice that the expression $\sin(\al/2)$ in \eqref{F17} was first
obtained by Szeg\H{o} \cite{Sz2}, where he calculated it as the Chebyshev constant of the arc $\G_{2\al}(\pi/2)$, then it was deduced
by Rosenblatt \cite{Ros}, as the capacity of $\G_{2\al}(\pi/2)$.}
\end{rem}
\begin{exa}
\label{ex2}
{\rm Let $\G_{2\al}(\al)$ be an arc of length $2\al$, defined by \eqref{arc}:
$\G_{2\al}(\al)=\{e^{i\theta}: \, \theta\in[0,2\al]\},$
and let $\G(2,\al)$ be the preimage of the arc $\G_{2\al}(\al)$ under the mapping
$p(z)=z^2$. It can be shown (see Babayan et al. \cite{BGT} for details)
that the set $\G(2,\al)$ is the union of two closed arcs of equal
lengths $\al$, symmetrically located with respect to the center of the unit circle (see Figure \ref{fig1}b):
\beq
\label{FF2}
\G(2,\al)=\{e^{i\om}: \, \om\in[-\pi,-\pi+\al]\cup [0,\al]\}.
\eeq
Then, by the Fekete theorem (see Proposition \ref{FT}) and formula \eqref{F17a},
for the transfinite diameter $\tau(\G(2,\al))$ we have
\beq
\nonumber
\tau(\G(2,\al)) =[\tau(\G_{2\al}(\al))]^{1/2}=\left(\sin(\al/2)\right)^{1/2}.
\eeq
The above result can easily be extended to the case of $k$ ($k>2$) arcs.
Let $\G(k,\al)$ be the union of $k$ ($k\in\mathbb{N},\, k\geq2$)
closed arcs of equal lengths $\al$, which are symmetrically located on the
unit circle (the arcs are assumed to be equidistant).
It can be shown that the set $\G(k,\al)$
is the preimage (to within rotation) under the mapping $p(z)=z^k$ of the arc $\G_{k\al}(k\al/2)$ of length $k\al$ defined by \eqref{arc}.
Therefore, by Fekete's formula \eqref{F0} 
and the invariance
property of the transfinite diameter with respect to rotation (see Proposition \ref{pp1}(b)),
for the transfinite diameter $\tau(\G(k,\al))$, we have
\beq
\label{F18}
\tau(\G(k,\al)) =\left(\sin(k\al/4)\right)^{1/k}.
\eeq}
\end{exa}
\begin{exa}
\label{ex3}
{\rm Let $\al>0,$ $\de\geq 0$ and $\al+\de\leq\pi$. Consider the set
\beq
\label{F188}
\G_{\al,\de}(\theta_0):=\G_{\al+\de}(\theta_0)\setminus\G_{\de}(\theta_0)
\eeq
consisting of the union of two arcs of the unit circle of lengths $\al$, the distance of
which (over the circle) is equal to $2\de$. Define (see Fig. 2a)):
\beq
\label{F18G}
\G_{\al,\de}:=\G_{\al,\de}(0)=\{e^{i\theta}: \, \theta\in
[-(\de+\al),-\de]\cup [\de,\de+\al]\}.
\eeq
Then the projection $\G_{\al,\de}^x$ of $\G_{\al,\de}$ onto the real axis
is the segment $\G_{\al,\de}^x=[\cos(\al+\de),\cos\de]$, and by Proposition \ref{pp1}(e)
for the transfinite diameter $\tau(\G_{\al,\de}^x)$ we have
$$\tau(\G_{\al,\de}^x) =\frac{\cos\de-\cos(\al+\de)}4
=\frac{\sin(\al/2)\sin(\al/2+\de)}2.$$
Hence, according to Robinson's formula \eqref{F10}, for the transfinite diameter $\tau(\G_{\al,\de})$,
we obtain
\beq
\label{F19}
\tau(\G_{\al,\de})= [2\tau(\G_{\al,\de}^x)]^{1/2} =\left(\sin(\al/2)\sin(\al/2+\de)\right)^{1/2}.
\eeq
In view of Proposition \ref{pp1}(b), from \eqref{F19} for any
$\theta_0\in[-\pi,\pi)$ we have
\beq
\label{F19a}
\tau(\G_{\al,\de}(\theta_0))=\left(\sin(\al/2)\sin(\al/2+\de)\right)^{1/2}.
\eeq}
\end{exa}
Observe that for $\de=0$ we have $\G_{\al,\de}(\theta_0)=\G_{2\al}(\theta_0)$ (see \eqref{arc} and \eqref{F188}), and formula \eqref{F19a} becomes \eqref{F17a}.
\begin{figure}[ht]%
\centering
\includegraphics[width=0.8\textwidth]{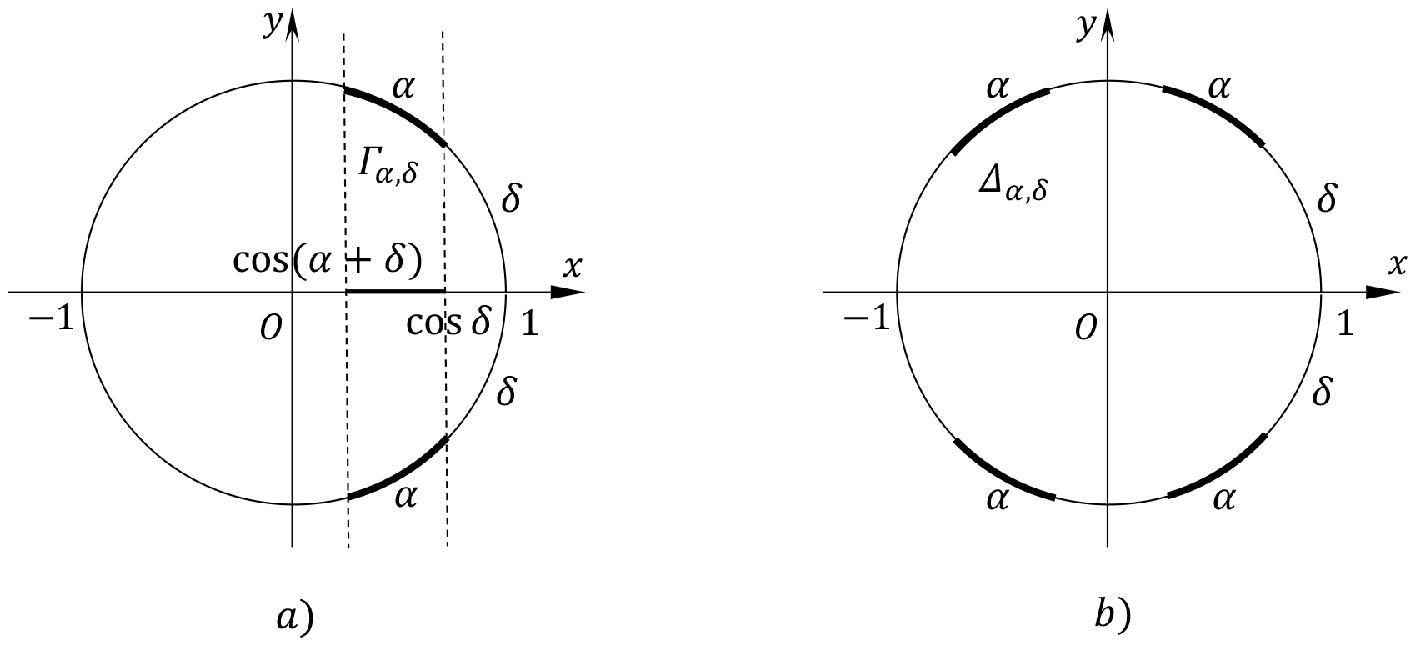}
\caption{a) The set $\G_{\al,\de}$. b) The set $\De_{\al,\de}$.}
\label{fig2}
\end{figure}

\begin{exa}
\label{ex4}
{\rm Let the arc $\G_{\al,\de}$ be as in Example \ref{ex3} (see \eqref{F18G})
with $\al,\de$
satisfying $\al+\de\leq\pi/2$, that is, $\G_{\al,\de}$ is a subset of the
right semicircle $\mathbb{T}$. Denote by $\G'_{\al,\de}$ the symmetric to
$\G_{\al,\de}$ set with respect to $y$-axis, that is,
$$\G'_{\al,\de}:=\{e^{i\theta}: \, \theta\in
[-\pi+\de,-\pi+(\de+\al)] \cup [\pi-(\de+\al),\pi-\de]\}.$$
Define
$\Delta_{\al,\de}:=\G_{\al,\de}\cup\G'_{\al,\de}$,
and observe that the set $\Delta_{\al,\de}$ consists of four arcs of equal lengths $\al$,
which are symmetrically located with respect to both axes (see Figure \ref{fig2}b)).
It is easy to see that the set
$\Delta_{\al,\de}$ is the preimage (to within rotation) of the set $\G_{2\al,2\de}$
under the mapping $p(z)=z^2$.
Hence, according to Fekete's formula \eqref{F0} and
\eqref{F19}, for the transfinite diameter
$\tau(\Delta_{\al,\de})$, we obtain
\beq
\label{F199}
\tau(\Delta_{\al,\de})= \left(\tau(\G_{2\al,2\de})\right)^{1/2}
=\left(\sin\al\sin(\al+2\de)\right)^{1/4}.
\eeq
Denote by $\Delta_{\al,\de}(\theta_0)$ the image of the set $\Delta_{\al,\de}$
under mapping $q(z)=e^{i\theta_0}z$, that is, under the rotation by the central angle
$\theta_0$ around the origin. Then, in view of Proposition \ref{pp1}(b)),
from \eqref{F199} for any $\theta_0\in[-\pi,\pi)$ we have
\beq
\label{F199a}
\tau(\Delta_{\al,\de}(\theta_0))=\left(\sin\al\sin(\al+2\de)\right)^{1/4}.
\eeq
}
\end{exa}

\ssn{A consequence of Theorem \ref{BB1}}

Now we apply Theorem \ref{BB1} to obtain
the asymptotic behavior of the prediction error $\si_n(f)$ in the
cases where the spectrum of a stationary process $X(t)$ is as in
Examples \ref{ex1}-\ref{ex4}.
Specifically, putting together Theorem \ref{BB1} and Examples \ref{ex1}-\ref{ex4},
we obtain the following result.

\begin{thm}[Babayan et al. \cite{BGT}]
\label{T5.10}
Let $\ol E_f$ be the support of the spectral density $f$ of a stationary process $X(t)$, and let $f>0$ a.e. on $\ol E_f$. Then for the prediction error $\sigma_n(f)$ the following assertions hold.
\begin{itemize}
\item[{\rm(a)}] If $\ol E_f= \G_{2\al}(\theta_0)$, where $\G_{2\al}(\theta_0)$
is as in Example \ref{ex1}, then
$$
\lim_{n\to\f}\sqrt[n]{\sigma_n(f)}  
=\sin(\al/2).$$
\item[{\rm(b)}] If $\ol E_f= \G(k,\al)$, where $\G(k,\al)$ is as in Example \ref{ex2}, then
$$
\lim_{n\to\f}\sqrt[n]{\sigma_n(f)}  
=\left(\sin(k\al/4)\right)^{1/k}.
$$
\item[{\rm(c)}] If $\ol E_f= \G_{\al,\de}(\theta_0)$, where $\G_{\al,\de}(\theta_0)$
is as in Example \ref{ex3}, then
$$
\lim_{n\to\f}\sqrt[n]{\sigma_n(f)} 
=\left(\sin(\al/2)\sin(\al/2+\de)\right)^{1/2}.
$$
\item[{\rm(d)}] If $\ol E_f= \Delta_{\al,\de}(\theta_0)$, where $\Delta_{\al,\de}(\theta_0)$,
is as in Example \ref{ex4}, then
$$
\lim_{n\to\f}\sqrt[n]{\sigma_n(f)}  
=\left(\sin\al\sin(\al+2\de)\right)^{1/4}.
$$
\end{itemize}
\end{thm}
\begin{rem}
{\rm The assertion (a) is a slight extension of the Rosenblatt relation \eqref{nd2}.
The assertion (c) is an extension of assertion (a), which reduces to assertion (a) if $\de=0$.}
\end{rem}

\sn{Davisson's theorem and its extension}
\label{dn1}

In this section, we consider a question of bounding the prediction
error $\si^2_n(f)$.
Using constructive methods, Davisson \cite{Dav1} obtained an upper bound
(rather than an asymptote) for the prediction
error $\si^2_n(f)$ without imposing continuity requirement on the
spectral density $f(\la)$. Specifically, in Davisson \cite{Dav1} was proved the following result:
\begin{thm}[Davisson \cite{Dav1}]
\label{D1}
Let the spectral density $f(\la)$, $\la\in[-\pi,\pi]$ of the process $X(t)$
be identically zero on a closed interval of length $2\pi-2\alpha$, $0<\alpha<\pi$.
Then for the prediction error $\si^2_n(f)$ the following inequality holds:
\beq
\label{ndd1}
\si^2_n(f) \le 4c\left(\sin(\alpha/2)\right)^{2n-2},
\eeq
where $c=r(0)$ and $r(\cdot)$ is the covariance function of $X(t)$
(see formula \eqref{mo1}).
\end{thm}
The theorem that follows, proved in Babayan and	Ginovyan \cite{BG-2},
extends Davisson's theorem to the case where the spectrum
of the process $X(t)$ consists of a union of two equal arcs.

Let $\al>0,$ $\de\geq 0$ and $\al+\de\leq\pi$, and let $\G_{\al,\de}$ be
the set defined by \eqref{F18G}.
Recall that $\G_{\al,\de}$
is the union of two arcs of the unit circle of lengths $\al$, the distance between which (over the circle) is equal to $2\de$ (see Example \ref{ex3}
and Figure \ref{fig2}a)).
\begin{thm}[Babayan and	Ginovyan \cite{BG-2}]
\label{D1E}
Let the spectral density $f(\la)$, $\la\in[-\pi,\pi]$ of the process $X(t)$ vanish outside the set $\G_{\al,\de}$.
Then for the prediction error $\si^2_n(f)$ the following inequality holds:
\beq
\label{d1e}
\si^2_n(f) \le 4c\left(\sin(\alpha/2)\right)^{n-1}\left(\sin(\alpha/2+\de)\right)^{n-1},
\eeq
where $c$ is as in Theorem \ref{D1}.
\end{thm}
\begin{rem}
\label{ed2}
{\rm For $\de=0$ the set $\G_{\al,\de}$ defined by \eqref{F18G}
is an arc of length $2\al$, and, in this case, the inequality \eqref{d1e}
becomes Davisson's inequality \eqref{ndd1}.}
\end{rem}

\sn{Extensions of Rosenblatt's second theorem}
\label{d2}

In this section, we analyze the asymptotic behavior of the prediction error
in the case where the spectral density $f(\la)$, $\la\in[-\pi,\pi]$
of the model is strictly positive
except one or several points at which it has a very high order
contact with zero so that the Szeg\H{o} condition \eqref{S} is violated.

Based on Rosenblatt's result for this case, namely Theorem \ref{R2},
we can expect that for any deterministic process with spectral density possessing
a singularity of the type \eqref{tt2}, the rate of the prediction error $\si^2_n(f)$
should be the same as in \eqref{nd6}.
However, the method applied in Rosenblatt \cite{Ros}
does not allow to prove this assertion. In Babayan and Ginovyan \cite{BG-1, BG}
and in Babayan et al. \cite{BGT}, using a different approach,
Rosenblatt's second theorem was extended to broader classes of spectral densities.
To state the corresponding results, we first examine the asymptotic behavior
as $n\to\f$ of the ratio:
$$
\frac{\sigma_n^2(fg)}{\sigma_n^2(f)},
$$
where $g$ is a non-negative function.

To clarify the approach, we first assume that $f$ is the spectral
density of a nondeterministic process, in which case the geometric mean $G(f)$
is positive (see \eqref{c013} and \eqref{a2}).
We can then write
\beq
\label{k3}
\lim_{n\to\infty}\frac{\sigma_n^2(fg)}{\sigma_n^2(f)}
=\frac{\sigma_\infty^2(fg)}{\sigma_\infty^2(f)}=\frac{2\pi G(fg)}{2\pi G(f)}
=\frac{G(f)G(g)}{G(f)}=G(g).
\eeq
It turns out that under some additional assumptions imposed on functions
$f$ and $g$, the asymptotic relation (\ref{k3}) remains also valid in the
case of deterministic processes, that is, when 
$G(f)=0$.

\ssn{Preliminaries}
\label{pre}

In what follows we consider the class of {\it deterministic} processes
possessing spectral densities $f$ for which the sequence of
prediction errors $\{\si_n(f)\}$ is weakly varying
(see Definition \ref{kd1}),
and denote by $\mathcal{F}$ the class of the corresponding spectral densities:
\beq
\label{k04}
\mathcal{F}:=\left\{f\in L^1(\Lambda): \, \, f\geq 0, \, \, G(f)=0, \,\, \lim_{n\to\f} \frac{\sigma_{n+1}(f)}{\sigma_{n}(f)} =1\right\}.
\eeq
\begin{rem}
\label{vvv2}
{\rm According to Rakhmanov's theorem (Theorem \ref{Rah}),
a sufficient condition for $f\in\mathcal{F}$ is that
$f>0$ almost everywhere on $\Lambda$ and $G(f)=0$.
Thus, the considered class $\mathcal{F}$ includes all deterministic
processes ($G(f)=0$) with almost everywhere positive spectral densities ($f>0$ a.e.).
On the other hand, according to Theorem 3.2 and Remark \ref{vvv},
the class $\mathcal{F}$ does not contain spectral densities, which vanish on an entire segment of $\Lambda$
(or on an arc of the unit circle $\mathbb{T}$).
Also, from Theorem \ref{cG} and Remark \ref{vvv} we infer that
a necessary condition for $f\in\mathcal{F}$ is that
the spectrum $E_f$ is dense in $\Lambda$.}
\end{rem}

\begin{den}
\label{kd33}
Let $\mathcal{F}$ be the class of spectral densities defined by \eqref{k04}.
For $f\in \mathcal{F}$ denote by $\mathcal{M}_f$ the class of 
nonnegative functions $g(\la)$ $(\la\in\Lambda)$ satisfying the conditions:
$G(g)>0$, $fg\in L^1(\Lambda)$, and
\beq
\label{e5.6} 
\lim_{n\to\f}\frac{\si^2_{n}(fg)}{\si^2_{n}(f)} =G(g),
\eeq
that is,
\beq
\label{e5.6a}
\mathcal{M}_f:=\left\{g\geq 0, \, \, G(g)>0, \,\, fg\in L^1(\Lambda),
\, \, \lim_{n\to\f}\frac{\si^2_{n}(fg)}{\si^2_{n}(f)} =G(g)\right\}.
\eeq
\end{den}

The next proposition shows that the class $\mathcal{F}$
is close under multiplication by functions from the class $\mathcal{M}_f$.
\begin{pp}[Babayan and Ginovyan \cite{BG}]
\label{ppg1}
If $f\in \mathcal{F}$ and $g\in \mathcal{M}_f$, then $fg\in \mathcal{F}$.
\end{pp}
The next result shows that the class $\mathcal{M}_f$ in a certain sense
is close under multiplication.
\begin{pp}[Babayan and Ginovyan \cite{BG}]
\label{ppg}
Let $f\in \mathcal{F}$. If $g_1\in \mathcal{M}_f$ and $g_2\in \mathcal{M}_{fg_1}$, then $g:=g_1g_2\in \mathcal{M}_f$ and $fg\in \mathcal{F}$.
In particular, if $g\in \mathcal{M}_f\cap\mathcal{M}_{fg}$,
then $g^2\in \mathcal{M}_f$.
\end{pp}

In the next definition we introduce certain classes of bounded functions.
\begin{den}
\label{kd2}
We define the class $B$ to be the set of all nonnegative, Riemann integrable on
$\Lambda=[-\pi,\pi]$ functions $h(\lambda)$. Also, we define the following subclasses:
\beq
\label{A}
B_+:= \{h \in B: \, h(\lambda)\geqslant m\},\q
B^-:= \{h \in B: \, h(\lambda)\leqslant M\}, \q B_+^-:=B_+\cap B^-,
\eeq
where $m$ and $M$ are some positive constants.
\end{den}

In the next proposition we list some obvious properties of the 
classes $B_+$, $B^-$ and $B_+^-$.
\begin{pp}
\label{p3.3}
The following assertions hold.
\begin{itemize}
\item[a)] If $h\in B_+ (B^-)$, then $1/h\in B^- (B_+)$.

\item[b)] If $h_1, h_2 \in B_+ (B^-)$, then $h_1+ h_2 \in B_+ (B^-)$ and
$h_1 h_2 \in B_+ (B^-)$.

\item[c)] If $h_1, h_2 \in B^-$ and $h_1/h_2$ is bounded,
then $h_1/h_2 \in B^-$.

\item[d)] If $h_1, h_2 \in B_+^-$, then $h_1+ h_2 \in B_+^-$,
 $h_1 h_2 \in B_+^-$ and $h_1/h_2 \in B_+^-$.
\end{itemize}
\end{pp}

In the next proposition we list some properties of weakly varying
sequences for functions from the above defined
classes $B_+$, $B^-$ and $B_+^-$ (see Babayan et al. \cite{BGT}).
\begin{pp}
\label{p3.3v}
Let the spectral density $f$ be such that the sequence $\sigma_n(f)$
is weakly varying. The following assertions hold.
\begin{itemize}
\item[a)] If $g\in B_+^-$, then the sequence $\sigma_n(fg)$ is also weakly varying.

\item[b)] If $g\in B^-$ with $G(g)=0$, then $\sigma_n(fg)=o(\sigma_n(f))$ as $n\to\infty$.
Thus, multiplying singular spectral densities we obtain a spectral density
with higher 'order of singularity'.

\item[c)] If $g\in B_+$ with $G(g)=\infty$, and $fg\in B$,
then $\sigma_n(f)=o(\sigma_n(fg))$ as $n\to\infty$.
\end{itemize}
\end{pp}

\ssn{Extensions of Rosenblatt's second theorem}
\label{fdR}
The following theorem, proved in Babayan et al. \cite{BGT},
describes the asymptotic behavior of the ratio $\sigma_n^2(fg)/\sigma_n^2(f)$
as $n\to\f$, and essentially states that if the spectral density $f$
is from the class $\mathcal{F}$ (see \eqref{k04}), and $g$
is a nonnegative function, which can have {\it polynomial} type
singularities, then the sequences $\{\si_n(fg)\}$ and
$\{\si_n(f)\}$  have the same asymptotic behavior as $n\to\f$ up to a
positive numerical factor.

\begin{thm}[Babayan et al. \cite{BGT}]
\label{sT1}
Let $f$ be an arbitrary function from the class $\mathcal{F}$, and let
$g$ be a function of the form:
\beq
\label{g}
g(\la)=h(\la)\cd\frac{t_1(\la)}{t_2(\la)}, \q \la\in\Lambda,
\eeq
where  $h\in B_+^-$,  $t_1$ and $t_2$ are nonnegative trigonometric polynomials, such that $fg\in L^1(\Lambda)$.
Then $g\in\mathcal{M}_f$ and $fg\in\mathcal{F}$, that is, $fg$ is the spectral density of a deterministic process with weakly varying prediction error,
and the relation \eqref{e5.6} holds.
\end{thm}

In view of Remark \ref{vvv2}, as a consequence of Theorem \ref{sT1}
we obtain the following result.
\begin{cor}
\label{c4.0}
Let the spectral density $f$ of a deterministic process $X(t)$ be a.e.
positive, and let $g$ be as in Theorem \ref{sT1}.
Then $g$ is the spectral density of a nondeterministic process and
the relation \eqref{e5.6} holds.
\end{cor}

As an immediate consequence of Theorem \ref{sT1} and Proposition \ref{p4.1}(d),
we have the following result.
\begin{cor}
\label{c4.1}
Let the functions $f$ and $g$ be as in Theorem \ref{sT1}.
Then the sequence $\sigma_n(fg)$ is also weakly varying.
\end{cor}

The theorems that follow extend the above stated Theorem \ref{sT1}
to a broader class of
spectral densities, for which the function $g$ can have
{\it arbitrary power type singularities}.

\begin{thm}[Babayan and Ginovyan \cite{BG}]
\label{T1}
Let $f$ be an arbitrary function from the class $\mathcal{F}$, and let
$g$ be a function of the form:
\beq
\label{e5.5}
g(\la)=h(\la)\cd |t(\la)|^\al, \q\al>0, \, \, \la\in\Lambda,
\eeq
where $h\in B_+^-$ and $t$ is an arbitrary trigonometric polynomial.
Then $g\in\mathcal{M}_f$ and $fg\in\mathcal{F}$, that is, $fg$ is the spectral density of a deterministic process with weakly varying prediction error,
and the relation \eqref{e5.6} holds.
\end{thm}

Using inductive arguments and Theorem \ref{T1} we can state the following result.
\begin{cor}
\label{ec5.2}
The conclusion of Theorem \ref{T1} remains valid if the function  $g$
has the following form:
\beq
\nonumber
g(\la)=h(\la)\cd |t_1(\la)|^{\al_1}\cd |t_2(\la)|^{\al_2}\cd\cdots\cd
|t_m(\la)|^{\al_m}, \q \la\in\Lambda,
\eeq
where  $h\in B_+^-$, \, $t_1, t_2, \ldots, t_m$ are arbitrary
trigonometric polynomials, $\al_1, \al_2, \ldots, \al_m$ are arbitrary
positive numbers, and $m\in\mathbb{N}$.
\end{cor}
\begin{thm}[Babayan and Ginovyan \cite{BG}]
\label{T2}
Let $f$ be an arbitrary function from the class $\mathcal{F}$, and let
$g$ be a function of the form:
\beq
\label{e5.5a}
g(\la)=h(\la)\cd t^{-\al}(\la), \q \al>0, \, \, \la\in\Lambda,
\eeq
where  $h\in B_+^-$ and $t$ is a nonnegative trigonometric polynomial.
Then the following assertions hold.
\begin{itemize}
	\item[(a)] $g\in\mathcal{M}_f$ and $fg\in\mathcal{F}$ provided that
	$\al\in\mathbb{Z}$ and $ft^{-\al}\in L^1(\Lambda)$.
	\item[(b)] $g\in\mathcal{M}_f$ and $fg\in\mathcal{F}$ provided that
	$\al\notin\mathbb{Z}$ and $ft^{-(k+1)}\in L^1(\Lambda)$, where $k:=[\al]$ is the integer part of $\al$.
\end{itemize}
\end{thm}

To state the next result we need the following definition.
\begin{den}
\label{kd3}
Let $E_1$ and $E_2$ be two numerical sets such that for any $x\in E_1$ and $y\in E_2$ we have $x<y$. We say that the sets $E_1$ and $E_2$ are separated from each other
if $\sup E_1<\inf E_2.$ Also, we say that a numerical set $E$ is separated from
infinity if it is bounded from above.
\end{den}
\begin{thm}[Babayan and Ginovyan \cite{BG}]
\label{T2c}
Let $f(\la)$ and $\hat f(\la)$ $(\la\in\Lambda)$ be spectral densities
of stationary processes satisfying the following conditions:
\begin{itemize}
\item[1)] $f, \hat f\in B^-$;

\item[2)] the functions  $f(\la)$ and $\hat f(\la)$ have $k$ common
essential zeros 
$$\la_1,\la_2,\ldots,\la_k\in\Lambda \q (-\pi<\la_1<\la_2<\cdots<\la_k\leq \pi, \, \, k\in\mathbb{N}),$$
that is,
\beq
\label{s3.7}
\lim_{\la\to\la_j}f(\la)=\lim_{\la\to\la_j}\hat f(\la)=0, \q j=1,2,\ldots,k;
\eeq
\item[3)] the functions  $f(\la)$ and $\hat f(\la)$ are infinitesimal of
the same order in a neighborhood of each point $\la_j$ $(j=1,2,\ldots,k)$, that is,
\beq
\label{s3.8}
\lim_{\la\to\la_j}\frac{\hat f(\la)}{f(\la)}=c_j>0, \q j=1,2,\ldots,k;
\eeq
\item[4)] the functions  $f(\la)$ and $\hat f(\la)$ are bounded away from zero
outside any neighborhood $O_\de(\la_j)$ $(j=1,2,\ldots,k)$, which is separated
from the neighboring zeros $\la_{j-1}$ and $\la_{j+1}$ of $\la_j$,
that is, there is a number $m:=m_\de>0$ such that $f(\la)\geq m$ and $\hat f(\la)\geq m$
for almost all $\la\notin\cup_{j=1}^kO_\de(\la_j)$.
Then the following assertions hold:
\end{itemize}
\begin{itemize}
\item[a)] $h(\la):=\frac{\hat f(\la)}{f(\la)}\in B_+^-;$
\item[b)] the processes with spectral densities $f$ and $\hat f$ either
both are deterministic or both are nondeterministic;
\item[c)] if one of the functions $f$ and $\hat f$ is from the class $\mathcal{F}$, then so is the other, and the following relation holds:
\beq
\label{g3.6e}
\lim_{n\to\f}\frac{\si_n^2(\hat f)}{\si_n^2(f)}=G(h)>0.
\eeq
\end{itemize}
\end{thm}

\begin{rem}
\label{r3.2}
{\rm The conditions of Theorem \ref{T2c} mean that the points $\la_j$
$(j=1,2,\ldots,k)$ are the only common zeros of functions $f(\la)$ and
$\hat f(\lambda)$. Besides, in the case of deterministic processes,
at least one of these zeros should be of sufficiently high order.
Also, notice that the conditions 1) and 4) of Theorem \ref{T2c} will be satisfied
if the functions $f(\la)$ and $\hat f(\lambda)$ are continuous on $\Lambda$.}
\end{rem}

\begin{thm}[Babayan and Ginovyan \cite{BG}]
\label{T1b}
Let $f$ be an arbitrary function from the class $\mathcal{F}$, and let
$g$ be a function of the form:
\beq
\label{e5.5b}
g(\la)=h(\la)\cd |q(\la)|^\al, \q \al\in\mathbb{R}, \, \, \la\in\Lambda,
\eeq
where  $h\in B_+^-$, $q$ is an arbitrary algebraic polynomial
with real coefficients, and $fg\in L^1(\Lambda)$.
Then $fg\in\mathcal{F}$ and $g\in\mathcal{M}_f$.
\end{thm}

Taking into account that the sequence $\{n^{-\al}, \, \, n\in\mathbb{N}, \, \al>0\}$
is weakly varying, as an immediate consequence of Theorems
\ref{T1}, \ref{T2}, \ref{T1b} or Corollary \ref{ec5.2},
we have the following result (see Babayan and Ginovyan \cite{BG}, and
Babayan et al. \cite{BGT}).
\begin{cor}
\label{c4.2}
Let the functions $f$ and $g$ satisfy the conditions of one of Theorems
\ref{T1}, \ref{T2}, \ref{T1b} or Corollary \ref{ec5.2},
and let $\sigma_n(f)\sim cn^{-\al}$ ($c>0, \al>0$) as $n\to\f$. Then
\beq
\label{s88}
\nonumber
\sigma_n(fg)\sim c G(g)n^{-\al} \q {\rm as} \q n\to\f,
\eeq
where $G(g)$ is the geometric mean of $g$.
\end{cor}

The next result, which immediately follows from Theorem \ref{R2} and Corollary \ref{c4.2},
extends Rosenblatt's second theorem (Theorem \ref{R2}) (see Babayan and Ginovyan \cite{BG}, and
Babayan et al. \cite{BGT}).
\begin{thm}[Babayan and Ginovyan \cite{BG}]
\label{sT2}
Let $f=f_ag$, where $f_a$ is defined by (\ref{nd4}),
and let $g$ be a function satisfying the conditions of one of Theorems
\ref{T1}, \ref{T2}, \ref{T1b} or Corollary \ref{ec5.2}. Then
\beq
\label{s888}
\nonumber
\de_n(f) = \si^2_n(f)\sim\frac{\Gamma^2\left(\frac{a+1}2\right)G(g)}
{\pi 2^{2-a}} \ n^{-a} 
\q {\rm as} \q n\to\f,
\eeq
where $G(g)$ is the geometric mean of $g$.
\end{thm}

We thus have the same limiting behavior for $\si^2_n(f)$ as in
the Rosenblatt's relation \eqref{nd6} up to an additional
positive factor $G(g)$.

\begin{rem}
\label{vvv3}
{\rm In view of Remark \ref{vvv2} it follows that all the above stated results
remain true if the condition $f\in\mathcal{F}$
is replaced by the following slightly strong but more constructive condition:
'{\it the spectral density $f$ is positive ($f>0$) almost everywhere on $\Lambda$ and $G(f)=0$'.}}
\end{rem}

\ssn{Examples}

In this section we discuss examples demonstrating the result stated in
Section \ref{fdR}.
In these examples we assume that $\{X(t),$ $t\in\mathbb{Z}\}$ is a stationary deterministic
process with a spectral density $f$ satisfying the conditions of Theorem \ref{sT1},
and the function $g$ is given by formula \eqref{g}.
To compute the geometric means we use the properties stated in Proposition \ref{p4.2}(a).
\begin{exa}
\label{ex41}
{\rm Let the function $g(\la)$ be as in \eqref{g} with $h(\la)=c>0$ and
$t_1(\la)=t_2(\la)=1$, that is, $g(\la)=c>0$. Then for the geometric mean
$G(g)$ we have
\beq
\label{g0}
G(g)=G(c)=c,
\eeq
and in view of \eqref{e5.6}, we get
\beq
\label{k71}
\nonumber
\lim_{n\to\f}\frac{\si^2_{n}(fg)}{\si^2_{n}(f)} =G(g)=c.
\eeq
Thus, multiplying the spectral density $f$ by a constant $c>0$
multiplies the prediction error by $c$.}
\end{exa}
\begin{exa}
\label{ex42}
{\rm Let the function $g$ be as in \eqref{g} with $h(\la)=e^{\I(\la)}$,
where $\I(\la)$ is an arbitrary odd function, and let
$t_1(\la)=t_2(\la)=1$, that is, $g(\la)=e^{\I(\la)}$.
Then for the geometric mean $G(g)$ we have
\beq
\label{g2}
G(g)=G(e^{\I(\la)})=\exp\left\{\frac1{2\pi}\inl\ln g(\la)\,d\la \right\}
=\exp\left\{\frac1{2\pi}\inl\I(\la)\,d\la \right\}=e^0=1,
\eeq
and in view of \eqref{e5.6}, we get
\beq
\label{k72}
\nonumber
\lim_{n\to\f}\frac{\si^2_{n}(fg)}{\si^2_{n}(f)} =G(g)=1.
\eeq
Thus, multiplying the spectral density $f$ by the function $e^{\I(\la)}$ with odd
$\I(\la)$ does not change the asymptotic behavior of the prediction error.}
\end{exa}

\begin{exa}
\label{ex43}
{\rm Let the function $g$ be as in \eqref{g} with $h(\la)=\la^2+1$ and
$t_1(\la)=t_2(\la)=1$, that is, $g(\la)=\la^2+1$. Then for the geometric mean $G(g)$
by direct calculation we obtain
\beq
\label{g3}
G(g)=\exp\left\{\frac1{2\pi}\inl\ln (\la^2+1)\,d\la \right\}
=\exp\{\ln(1+\pi^2)-2+\frac2\pi\arctan\pi\}\approx 3.3,
\eeq
and in view of \eqref{e5.6}, we get
\beq
\nonumber
\lim_{n\to\f}\frac{\si^2_{n}(fg)}{\si^2_{n}(f)} =G(g)
=\exp\{\ln(1+\pi^2)-2+\frac2\pi\arctan\pi\}\approx 3.3.
\eeq
Thus, multiplying the spectral density $f$ by the function $\la^2+1$
multiplies the prediction error approximately by 3.3.}
\end{exa}

\begin{exa}
\label{ex44}
{\rm Let the function $g$ be as in \eqref{g} with $h(\la)=t_2(\la)=1$,
and $t_1(\la)=\sin^{2k}(\la-\la_0)$, where $k\in\mathbb{N}$ and $\la_0$
is an arbitrary point from $[-\pi,\pi]$, that is, $g(\la)=\sin^{2k}(\la-\la_0)$.
To compute the geometric mean $G(g)$, we first find the algebraic
polynomial $s_2(z)$ in the Fej\'er-Riesz representation \eqref{ss2}
of the non-negative trigonometric polynomial $\sin^2(\la-\la_0)$ of degree 2.
For any $\la_0\in[-\pi,\pi]$ we have
\bea
\nonumber
\sin^2(\la-\la_0) =|\sin(\la-\la_0)|^2=
\left|\frac12(e^{2i(\la-\la_0)}-1)\right|^2
=\left|s_2(e^{i\la})\right|^2,
\eea
where
\beq
\label{k741}
s_2(z)=\frac12(e^{-2i\la_0}z^2-1).
\eeq
Therefore, by Proposition \ref{p4.2}(d) and \eqref{k741}, we have
\beq
\label{k742gg}
G(\sin^{2}(\la-\la_0))=|s_2(0)|^2=\left(\frac12\right)^2=\frac14.
\eeq
Now, in view of Proposition \ref{p4.2}(a) and \eqref{k742gg}, for the geometric mean of $g(\la)=t_1(\la)=\sin^{2k}(\la-\la_0)$ ($k\in\mathbb{N}$), we obtain
\beq
\label{g4}
G(g)=G(\sin^{2k}(\la-\la_0))=G^k(\sin^2(\la-\la_0))={4^{-k}},
\eeq
and in view of \eqref{e5.6}, we get
\beq
\label{k76}
\nonumber
\lim_{n\to\f}\frac{\si^2_{n}(fg)}{\si^2_{n}(f)} =G(g)=\frac1{4^k}.
\eeq
Thus, multiplying the spectral density $f$ by the non-negative trigonometric polynomial
$\sin^{2k}(\la-\la_0)$  of degree $2k$ ($k\in\mathbb{N}$), yields a
$4^k$-fold asymptotic reduction of the prediction error.}
\end{exa}
\begin{exa}
\label{ex45}
{\rm Let the function $g$ be as in \eqref{g} with $h(\la)=t_1(\la)=1$,
and $t_2(\la)=\sin^{2l}(\la-\la_0)$, where $l\in\mathbb{N}$ and $\la_0$
is an arbitrary point from $[-\pi,\pi]$, that is, $g(\la)=\sin^{-2l}(\la-\la_0)$.
Then, in view of the third equality in \eqref{gt1} and \eqref{g4}
for the geometric mean $G(g)$ we have
\beq
\label{g5}
G(g)=G(\sin^{-2l}(\la-\la_0))=G^{-1}(\sin^{2l}(\la-\la_0)) ={4^l},
\eeq
and in view of \eqref{e5.6}, we get
\beq
\label{k76g}
\nonumber
\lim_{n\to\f}\frac{\si^2_{n}(fg)}{\si^2_{n}(f)} =G(g)={4^l}.
\eeq
Thus, dividing the spectral density $f$ by the non-negative trigonometric polynomial
$\sin^{2l}(\la-\la_0)$  of degree $2l$ ($l\in\mathbb{N}$), yields
a $4^l$-fold asymptotic increase of the prediction error.}
\end{exa}
Notice that the values of the geometric mean $G(g)$ obtained in
\eqref{g4} and \eqref{g5} do not depend on the choice of the point
$\la_0\in[-\pi,\pi]$.

Putting together Examples \ref{ex41} - \ref{ex45} and using Proposition \ref{p4.2}(a)
we have the following summary example.
\begin{exa}
\label{ex46}
{\rm Let $\{X(t), \,t\in\mathbb{Z}\}$ be a stationary deterministic
process with a spectral density $f$ satisfying the conditions
of Theorem \ref{sT1}. Let $h(\la)=ce^{\I(\la)}(\la^2+1)$, $t_1(\la)=\sin^{2k}(\la-\la_1)$
and $t_2(\la)=\sin^{2l}(\la-\la_2)$, where $c$ is an arbitrary positive constant,
$\I(\la)$ is an arbitrary odd function and $\la_1, \la_2$ are arbitrary points from $[-\pi,\pi]$.
Let the function $g$ be defined as in \eqref{g}, that is,
\beq
\label{g6}
g(\la)=h(\la)\cd\frac{t_1(\la)}{t_2(\la)}=ce^{\I(\la)}(\la^2+1)
\frac{\sin^{2k}(\la-\la_1)}{\sin^{2l}(\la-\la_2)}.
\eeq
Then, in view of Proposition \ref{p4.2}(a) and relations \eqref{g0}--\eqref{g3}
and \eqref{g4}--\eqref{g6}, we have
\bea
\label{g7}
\nonumber
G(g)&=&
G(h)\frac{G(t_1)}{G(t_2)}
=G(c)G(e^{\I})G(\la^2+1)G(\sin^{2k}(\la-\la_1))G(\sin^{-2l}(\la-\la_2))\\
&=&(c)(1)\exp\{\ln(1+\pi^2)-2+\frac2\pi\arctan\pi\}(4^{-k})(4^{l})
\approx 3.3c4^{l-k},
\eea
and in view of \eqref{e5.6} and \eqref{g7}, we get
\beq
\nonumber
\lim_{n\to\f}\frac{\si^2_{n}(fg)}{\si^2_{n}(f)} =G(g)
\approx 3.3c4^{l-k}.
\eeq
}
\end{exa}

\begin{exa}
\label{sex1}
{\rm Let the function $g(\la)$ ($\la\in\Lambda$) be as in \eqref{e5.5} with $h(\la)=1$
and $t(\la)=\sin(\la-\la_0)$, where $\la_0$ is an arbitrary point from
$[-\pi,\pi]$, that is, $g(\la)=|\sin(\la-\la_0)|^{\al}$, $\al>0$.
Then, according to Example \ref{ex44}, 
for the geometric mean of $\sin^{2}(\la-\la_0)$ we have
\beq
\label{k742g}
G(\sin^{2}(\la-\la_0))=\frac14.
\eeq
According to Proposition \ref{p4.2}(a) and \eqref{k742g},
for the geometric mean of $g(\la)$, we obtain
\beq
\label{gg4}
G(g)=G(|\sin(\la-\la_0)|^{\al})=G\left(\left(\sin^{2}(\la-\la_0)\right)^{\al/2}\right)=
G^{\al/2}(\sin^2(\la-\la_0))=\frac1{2^\al},
\eeq
and in view of \eqref{e5.6}, we get
\beq
\nonumber
\lim_{n\to\f}\frac{\si^2_{n}(fg)}{\si^2_{n}(f)} =G(g)=\frac1{2^\al}.
\eeq
Thus, multiplying the spectral density $f(\la)$ by the function
$g(\la)=|\sin(\la-\la_0)|^{\al}$
yields a $2^\al$-fold asymptotic reduction of the prediction error.}
\end{exa}

\begin{exa}
\label{dex2}
{\rm Let the function $g(\la)$ be as in \eqref{e5.5b} with $h(\la)=1$
and $q(\la)=\la$, that is, $g(\la)=|\la|^\al$, $\al\in\mathbb{R}$.
By direct calculation we obtain
\beaa
\ln G(g)=\frac1{2\pi}\inl\ln |\la|^\al\,d\la =\frac\al\pi\int_0^\pi\ln\la\,d\la =\al\ln(\pi/e).
\eeaa

Therefore
\beaa
G(g)= \left(\pi/ e\right)^\al \approx (1.156)^\al,
\eeaa
and in view of \eqref{e5.6}, we get
\beq
\nonumber
\lim_{n\to\f}\frac{\si^2_{n}(fg)}{\si^2_{n}(f)} =G(g)
=\left(\frac\pi e\right)^\al \approx (1.156)^\al.
\eeq
Thus, multiplying the spectral density $f(\la)$ by the function
$g(\la)=|\la|^\al$ multiplies the prediction error asymptotically
by $(\pi/e)^\al \approx (1.156)^\al$.

It follows from Proposition \ref{pp3}(d) that the same asymptotic
is true for the prediction error with spectral density
$\bar g(\la)=|\la-\la_0|^\al$, $\la_0\in [-\pi,\pi]$.}
\end{exa}

\ssn{Rosenblatt's second theorem revisited}
\label{dex3}

\begin{figure}[ht]%
\centering
\includegraphics[width=0.7\textwidth]{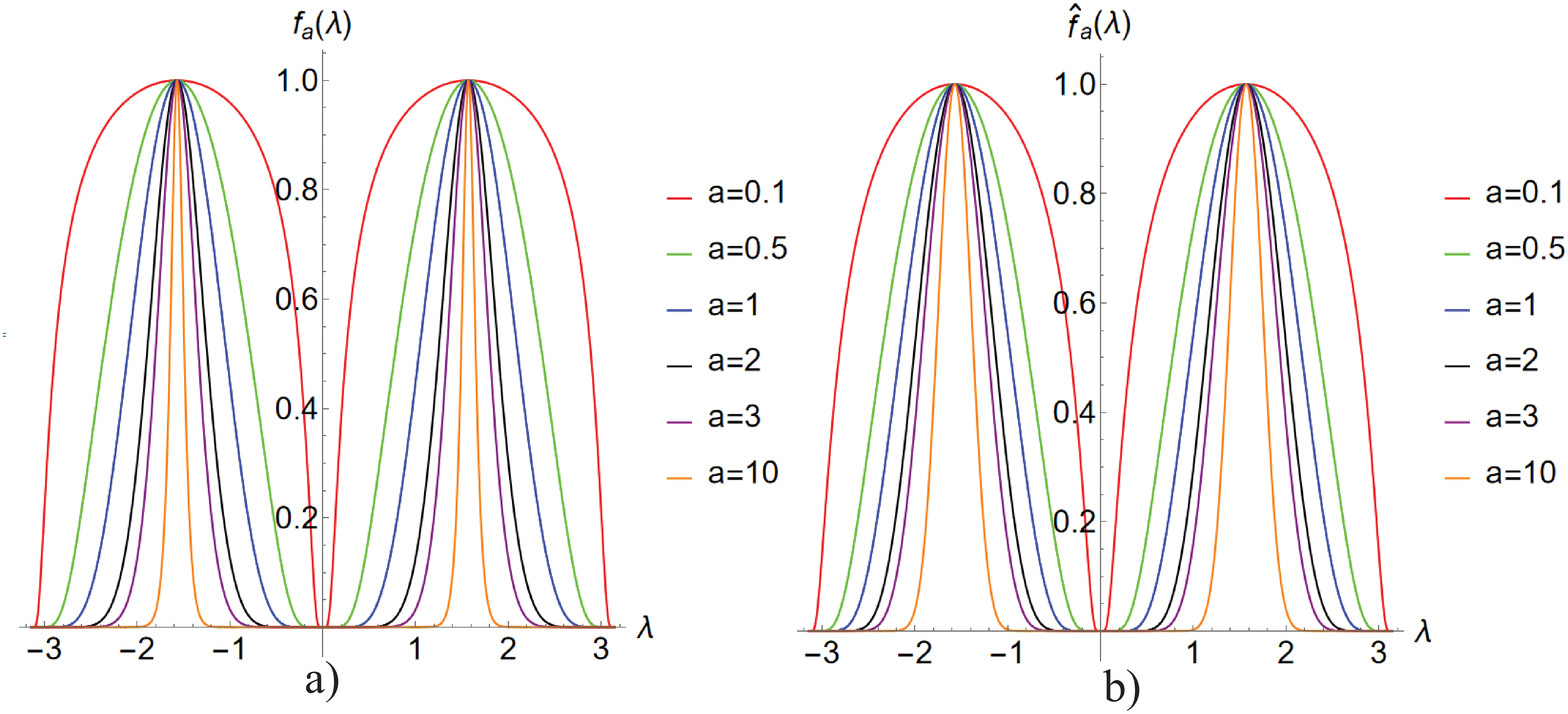}
\caption{a) Graph of the function $f_a(\lambda)$.
b) Graph of the function $\hat f_a(\lambda)$.}
\label{fig3}
\end{figure}

We first analyze the Pollaczek-Szeg\H{o} function
$f_a(\la)$ given by \eqref{nd4} (cf.  Pollaczek \cite{Pol} and Szeg\H{o} \cite{S1}). We have
\beq
\label{nd41}
f_a(\la)=\frac{2e^{2\la\varphi(\la)}e^{-\pi\varphi(\la)}}{e^{\pi\varphi(\la)}+e^{-\pi\varphi(\la)}}
=\frac{2e^{2\la\varphi(\la)}}{e^{2\pi\varphi(\la)}+1},
\q  0\leq\la\leq\pi, \q \varphi(\la):=\varphi_a(\la)=(a/2)\cot\la. 
\eeq
Observe that $\varphi(\la)\to+\f$ as $\la\to 0^+$, and we have 
\beq
\label{nd42}
\varphi(\la)\sim a/(2\la), \q e^{2\la\varphi(\la)}\sim e^a, \q
e^{2\pi\varphi(\la)}+1 \sim e^{a\pi/\la} \q {\rm as} \q \la\to0^+.
\eeq
Taking into account that $f_a(\la)$ is an even function, from \eqref{nd41}
and \eqref{nd42} we obtain the following asymptotic relation for $f_a(\la)$
in a vicinity of the point $\la=0$.
\beq
\label{nd43}
f_a(\la)\sim 2e^a\exp\left\{-{a\pi}/{|\la|}\right\}\q {\rm as}
\q \la\to0.
\eeq
Next, observe that $\varphi(\la)\to-\f$ as $\la\to \pi$, and we have
\beq
\label{nd44}
\varphi(\la)=-\varphi(\pi-\la) \sim (-a/2)(\pi-\la), \q 2\la\varphi(\la)\sim -a\pi/(\pi-\la), \q
{\rm as} \q \la\to\pi.
\eeq
In view of \eqref{nd41} and \eqref{nd44} we obtain the following asymptotic
of the function $f_a(\la)$ in a vicinity of the point $\la=\pi$.
\beq
\label{nd45}
f_a(\la)\sim 2e^{2\la\varphi(\la)} \sim 2\exp\left\{-{a\pi}/{(\pi-\la)}\right\}\q {\rm as}
\q \la\to\pi.
\eeq
Putting together \eqref{nd43} and \eqref{nd45}, and taking into account
evenness of $f_a(\la)$, we conclude that
\beq \label{t2}
f_a(\la)\sim
\left \{
\begin{array}{ll}
 2e^a\exp\left\{-{a\pi}/{|\la|}\right\} & \mbox{as $\la\to0$},\\
2\exp\left\{-{a\pi}/{(\pi-|\la|)}\right\} & \mbox{as $\la\to\pm\pi$},
\end{array}
\right.
\eeq
Thus, the function $f_a(\la)$ is positive everywhere except for points
$\la=0, \pm\pi,$ and has a very high order of contact with zero at these points, so that Szeg\H{o}'s condition \eqref{S} is violated implying that $G(f_a)=0$. Also, observe that $f_a(\la)$ is infinitely differentiable at
all points of the segment $[-\pi,\pi]$ including the points $\la=0, \pm\pi,$
and attains it maximum value of 1 at the points $\pm\pi/2$.
For some specific values of the parameter $a$ the graph of the function $f_a(\la)$ is represented in Figure \ref{fig3}a).

For $a>0$ and  $\la\in [-\pi,\pi]$, consider the pair of functions
$\hat f_1(\la)$ and $\hat f_2(\la)$ defined by formulas:
\beq
\label{ff}
\hat f_1(\la):=\exp\left\{-{a\pi}/{|\la|}\right\}, \q
\hat f_2(\la):=\exp\left\{-{a\pi}/{(\pi-|\la|)}\right\}.
\eeq

\begin{figure}[ht]%
\centering
\includegraphics[width=0.7\textwidth]{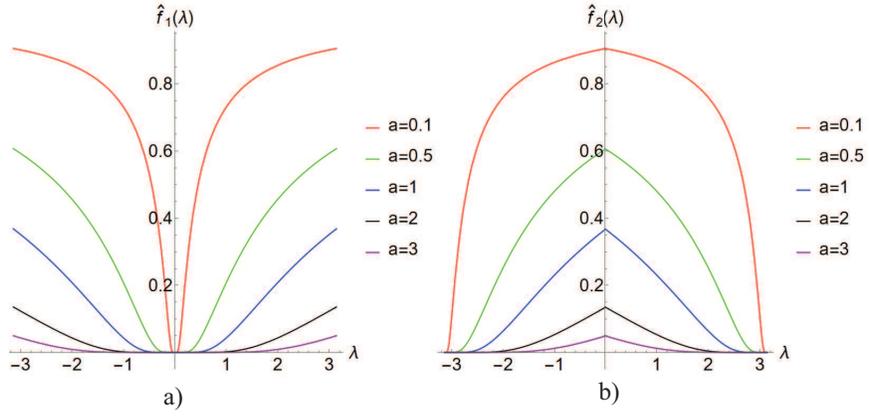}
\caption{a) Graph of the function $\hat f_1(\lambda)$.
b) Graph of the function $\hat f_2(\lambda)$.}
\label{fig4}
\end{figure}

Observe that the function $\hat f_1(\la)$ is positive everywhere except
for point $\la=0$ at which it has the same order of contact with zero
as $f_a(\la)$, and hence $G(\hat f_1)=0$. Also, $\hat f_1(\la)$ is infinitely differentiable at
all points of the segment $[-\pi,\pi]$  except for the points
$\la=\pm\pi,$ where it attains its maximum value equal to $e^{-a}$.
As for the function $\hat f_2(\la)$, it is positive everywhere except
for points $\la=\pm\pi,$ at which it has the same order of contact with
zero as $f_a(\la)$,
and hence $G(\hat f_2)=0$. Also, $\hat f_2(\la)$ is infinitely differentiable
at all points of the segment $[-\pi,\pi]$  except for the point $\la=0$,
where it attains its maximum value equal to $e^{-a}$.
For some specific values of the parameter $a$ the graphs of functions
$\hat f_1(\la)$ and $\hat f_2(\la)$ are represented in Figure \ref{fig4}.

Denote by $\hat f_a(\la)$ the product of functions $\hat f_1(\la)$
and $\hat f_2(\la)$ defined in \eqref{ff} and normalized by the factor $e^{4a}$:
\beq
\label{snd5}
\hat f_a(\la):=e^{4a}\hat f_1(\la)\hat f_2(\la)
=e^{4a}\exp\left\{-{a\pi^2}/{(|\la|(\pi-|\la|))}\right\},
\eeq
and observe that $\hat f_a(\la)$ behaves similar to $f_a(\la)$.
Indeed, the function $\hat f_a(\la)$ also is positive everywhere except for points $\la=0, \pm\pi,$ it is infinitely differentiable at all points of the segment $[-\pi,\pi]$ including the points $\la=0, \pm\pi,$ and attains it maximum value of 1 at the points $\pm\pi/2$.
Also, in view of \eqref{t2} and \eqref{snd5}, at points $\la=0, \pm\pi$ the function $\hat f_a(\la)$ has the same order of zeros as $f_a(\la)$, and hence $G(\hat f_a)=0$.
Thus, the process $X(t)$ with spectral density $\hat f_a(\la)$ is deterministic.
For some specific values of the parameter $a$ the graph of the function $\hat f_a(\la)$ is represented in Figure \ref{fig3}b).

The functions $f_a(\la)$ and $\hat f_a(\la)$ defined
by \eqref{nd4} and \eqref{snd5}, respectively, satisfy the conditions of
Theorem \ref{T2c}.
Therefore, we have (see \eqref{g3.6e})
\beq
\label{snd6}
\lim_{n\to\f}\frac{\si_n^2(\hat f_a)}{\si_n^2(f_a)}=G(\hat f_a/f_a):=\hat C(a)>0.
\eeq
In view of \eqref{nd6} and \eqref{snd6} we have
\beq
\label{snd7}
\si^2_n(\hat f_a)\sim C(a) \cd n^{-a}
\q {\rm as} \q n\to\f.
\eeq
where
\beq
\label{snd8}
C(a):=\frac{\Gamma^2\left(({a+1)}/2\right)\hat C(a)}{\pi 2^{2-a}}. 
\eeq

The values of the constants $\hat C(a)$ and $C(a)$ for some
specific values of the parameter $a$ are given in Table 1.

\begin{table}
\centering
\caption{The values of constants $\hat C(a)$ and $C(a)$}
\label{tab1}
\begin{tabular}{llll}
\toprule
{$a$} & {$\frac{\Gamma^2\left((a+1)/2\right)}{\pi 2^{2-a}}$} & {$\hat C(a)$} & {$C(a)$}\\
\midrule
0.1  & 0.223  & 0.797   & 0.178  \\
0.5  & 0.169  & 1.113   & 0.188  \\
    1.0  & 0.159  & 2.545   & 0.406  \\
    1.5  & 0.185  & 6.446   & 1.193  \\ 
    2.0  & 0.250  & 16.830  & 4.214   \\
    3.0  & 0.637  & 119.220 & 76.379  \\
    3.3  & 0.902   & 215.715 & 194.656  \\ \midrule
    3.4  & 1.020   & 263.173 & 268.375   \\
    5.0  & 10.186  & 6128.990& 62429.000  \\
    10.0   & 223256 & 1.104 $\cdot 10^8$ & 2.428 $\cdot 10^{13}$\\
\bottomrule
\end{tabular}
\end{table}

Now we compare the prediction errors $\si^2_n(\hat f_1)$ and
$\si^2_n(\hat f_2)$ with $\si^2_n(f_a)$.
To this end, observe first that the function
$g_1(\la):=f_a(\la)/\hat f_1(\la)$ has a very high order of contact
with zero at points $\la=\pm\pi,$ so that Szeg\H{o}'s condition \eqref{S}
is violated implying that $G(g_1)=0$. Besides, the function $g_1(\la)$ is continuous
on $[-\pi,\pi]$, and hence $g_1\in B^-$. Therefore, according to
Proposition \ref{p3.3v} b), we have
\beq
\label{ff1}
\si^2_n( f_a)= o\left(\si^2_n(\hat f_1)\right)
\q {\rm as} \q n\to\f.
\eeq
Similar arguments applied to the function $f_2(\la)$ yield
\beq
\label{ff2}
\si^2_n( f_a)= o\left(\si^2_n(\hat f_2)\right)
\q {\rm as} \q n\to\f.
\eeq

The relations  \eqref{ff1} and \eqref{ff2} show that the rate of
convergence to zero of the prediction errors $\si^2_n(\hat f_1)$
and $\si^2_n(\hat f_2)$ is less than the one for $\si^2_n(f_a)$,
that is, the power rate of convergence $n^{-a}$ (see \eqref{snd7}).
Thus, the rate of convergence $n^{-a}$ is due to the joint contribution
of all zeros $\la=0,\pm\pi$ of the function $f_a(\la)$,
whereas
each of these zeros separately does not guarantee the rate of convergence $n^{-a}$.

\section{An Application. Asymptotic behavior of the extreme
eigenvalues of truncated Toeplitz matrices}
\label{Ap}
In this section we analyze the relationship between the rate of
convergence to zero of the prediction error $\si_{n}^2(f)$ and the
minimal eigenvalue of a truncated Toeplitz matrix generated by the
spectral density $f$, by showing how it is possible
to obtain information in both directions.

The problem of asymptotic behavior of the extreme eigenvalues of truncated
(finite sections) Toeplitz matrices goes back to the classical works by
Kac, Murdoch and Szeg\H{o} \cite{KMS}, Parter \cite{Part}, Widom \cite{Wid},
and Chan \cite{Ch},
where the problem was studied for truncated Toeplitz matrices generated by continuous and continuously differentiable
functions (symbols). Since then the problem for various classes of symbols
was studied by many authors. For instance, Pourahmadi \cite{Po2},
Serra \cite{Se1,Se2}, and Babayan and Ginovyan \cite{BG-2} considered the problem in the case where the symbol of Toeplitz matrix is not
(necessarily) continuous nor differentiable (see also B\"ottcher and Grudsky \cite{BoG}).
In this section, we review and summarize some known results from the
above cited references and state some new results.

\subsection{Extreme eigenvalues of truncated Toeplitz matrices}
\label{Ap1}

Let $f(\la)$ be a real-valued Lebesgue integrable function defined on $\La:=[-\pi,\pi]$,
$T_n(f):=||r_{k-j}||,$ ${j,k=0,1,...,n}$, be the truncated Toeplitz matrix generated by the
Fourier coefficients of $f$, and let $\la_{1,n}(f) \le \la_{2,n}(f) \le \cdots\la_{n+1,n}(f)$ be the eigenvalues of $T_n(f)$.
We denote by $m_f:={\rm ess \ inf} f$ and $M_f:={\rm ess \ sup} f$ the essential minimum and the
essential maximum of $f$, respectively.
We will refer to $f(\la)$ as a {\it symbol} for the Toeplitz matrix $T_n(f)$.

We first recall Szeg\H{o}'s distribution theorem (see, e.g., Grenander and Szeg\H{o} \cite{GS}, p. 64-65).
\begin{thm}
\label{SDT}
For every continuous function $F$ defined in $[m_f, M_f]$
the following asymptotic relation holds:
\beq
\label {ea1}
\lim_{n\to\infty}\frac1{n+1}\sum_{k=1}^{n+1}F(\la_{k,n}(f))
= \frac1{2\pi}\inl F(f(u))\,du.
\eeq
Moreover, the spectrum of $T_n(f)$ is contained in $(m_f, M_f)$, and
\beq
\label {ea2}
\lim_{n\to\infty}\la_{1,n}(f) = m_f \q {\rm and} \q \lim_{n\to\infty}\la_{n+1,n}(f) = M_f.
\eeq
\end{thm}
The problem of interest is to describe the rate of convergence in \eqref{ea2},
depending on the properties of the symbol $f$.
In the following, without loos of generality, we assume that $m_f:={\rm ess \ inf} f=0$.
Also, we study the asymptotic behavior of the minimum eigenvalue
of $T_n(f)$, for the maximum eigenvalue it is sufficient to consider the
minimum eigenvalue of the matrix $T_n(-f)$.

The rate of convergence of extreme eigenvalues
has been studied by Kac, Murdoch and Szeg\H{o} \cite{KMS}, Parter \cite{Part},
Widom \cite{Wid} and Chan \cite{Ch}, under the following regularity condition on $f$
(see, e.g., Grenander and Szeg\H{o} \cite{GS}, Section 5.4(a), p. 72).

\vskip2mm
\noindent
{\bf Condition A}. 
Let $f(\la)$ be real, continuous and periodic with period $2\pi$.
Let $\min f(\la)=f(\la_0)=0$ and let $\la=\la_0$ be the only value of $\la$
(mod $2\pi)$ for which
this minimum is attained. Moreover, let $f(\la)$ have continuous derivatives of order
$2k$ $(k\in \mathbb{N})$ in some neighborhood of $\la=\la_0$ with $f^{(2k)}(\la_0)\neq 0$.

\begin{thm}[Kac, Murdoch and Szeg\H{o} \cite{KMS}]
\label{KMS}
Under Condition A the following asymptotic relation holds:
\beq
\label{ea3}
\lambda_{1,n}(f) \simeq n^{-2k} \q {\rm as} \q n\to\infty.
\eeq
\end{thm}
Observe that Condition A is too restrictive which may be hard to
verify or may even be unsatisfied in some areas of application such
as prediction theory of stationary processes and signal processing,
where $f$ is viewed as a spectral density of a stationary process.
Another limitation is the possibility that $f$ has more than one zero.

By using new linear algebra tools, Serra \cite{Se1} has extended
Theorem \ref{KMS}, by proving that the rate of convergence of
$\lambda_{1,n}(f)$ depends only on the order of the zero of $f$,
but not (necessarily) on the smoothness of $f$ as it is required
in Theorem \ref{KMS}.
Moreover, in Serra \cite{Se3} Theorem \ref{KMS} was further extended
to the case of a function $f \in L^1[-\pi, \pi]$ having several
global minima by showing that the maximal order of the zeros of
function $f$ is the only parameter which characterizes the rate of
convergence of $\lambda_{1,n}(f)$.
In particular, in Serra \cite{Se1} was proved the following result.

\begin{pp}[Serra \cite{Se1}]
\label{SeT1}
Let $f$ be a nonnegative integrable function on $[-\pi,\pi]$.
If $f$ has a zero of order $2k$ at a point $\la=\la_0$, that is,
$f(\la)\asymp(\la-\la_0)^{2k}$, then $\lambda_{1,n}(f)\asymp n^{-2k}$.
\end{pp}

\subsection{The relationship between the prediction error and the minimal eigenvalue.}
\label{Ap2}

Let $X(t),$ $t = 0,\pm1,\ldots,$ be a stationary sequence
possessing a spectral density function $f(\la),$ $\la\in [-\pi, \pi],$
and let $\sigma_n^2(f)$ be the prediction error in predicting $X(0)$
by the past of $X(t)$ of length $n$ (see formula \eqref{OS1}).

The next proposition provides a relationship between the minimal eigenvalue
$\lambda_{1,n}(f)$ of a truncated Toeplitz matrix $T_n(f)$ generated
by spectral density $f$ and the prediction error $\sigma_n^2(f)$
(see Pourahmadi \cite{Po2} and Serra \cite{Se1}).

\begin{pp}
\label{lemma:1}
Let $f$, $\lambda_{1,n}(f)$ and  $\sigma_n^2(f)$ be as above.
Then for any $n\in \mathbb{N}$ the following inequalities hold:
\beq
\label{ea7}
\lambda_{1,n}(f) \leq \si_n^2(f) \leq M_f\frac{\lambda_{1,n}(f)}{\lambda_{1,n-1}(f)}.
\eeq
\end{pp}

The first inequality in \eqref{ea7} was proved in Pourahmadi \cite{Po2},
while the proof of the second inequality in \eqref{ea7} can be found in
Serra \cite{Se1}.

Recall that for a stationary process $X(t)$ with spectral density $f(\la)$
by $E_f$ we denote the spectrum of $X(t)$, that is,
$E_f:=\{\la: \, f(\la)>0\}$ (see \eqref{Sp}). Thus, the closure $\ol E_f$
of $E_f$ is the support of the spectral density $f$.
Also, by $\tau^*(E_f)$ we denote the outer transfinite diameter of the
set $E_f$ (defined by \eqref{tdo1}).

The following theorem is an immediate consequence of Theorem \ref{B3}(b)
and Proposition \ref{lemma:1} (cf. Pourahmadi \cite{Po2}.)

\begin{thm}\label{thm:1}
Let $f$, $\lambda_{1,n}(f)$ and $\tau^*(E_f)$ be as above. Then the following inequality holds:
\bea
\label{ea13}
\limsup_{n\to\f}\sqrt[2n]{\lambda_{1,n}(f)} \leq \tau^*(E_f).
\eea
\end{thm}

Thus, in order that the minimal eigenvalue $\lambda_{1,n}(f)$ should decrease
to zero at least exponentially as $n\to\infty$, it is sufficient that
the outer transfinite diameter of the spectrum of the process $X(t)$
be less than 1.
As such the continuity and differentiability of spectral density
$f$ are not required for the exponential rate of convergence of
the minimal eigenvalue $\lambda_{1,n}(f)$ to zero.

Now we proceed to discuss two specific models of deterministic processes, for which can be obtained more information on the rate of convergence of
the minimal eigenvalue $\lambda_{1,n}(f)$ to zero as $n\to\infty$
from that of prediction error $\si_n^2(f)$.
Notice that, for the first model, the spectral density $f$ of the
process is discontinuous and has uncountably many zeros, while, for the
second model, the function $f$ has a zero of exponential order at
points $0, \pm\pi$.
Therefore, in both cases, Condition A of Kac, Murdoch and Szeg\H{o}
is violated.

\ssn{A model with spectral density $f$ which is discontinuous,
zero on an interval, and positive elsewhere.}
Let $X(t)$ be a stationary process for which the support $\ol E_f$
of the spectral density $f(\la)$ is as in Examples \ref{ex1} - \ref{ex4}.
Then, we can apply Theorems \ref{T5.10} and \ref{thm:1}
to obtain asymptotic estimates for the the minimal eigenvalue $\lambda_{1,n}(f)$.

In the next theorem we state the result in the cases where the support
$\ol E_f$ of $f$ is as in Examples \ref{ex1} and \ref{ex3},
similar estimates can be stated in the cases where
$\ol E_f$ is as in Examples \ref{ex2} and \ref{ex4}.

\begin{thm}[Babayan and Ginovyan \cite{BG-2}]
\label{T6.5}
Let $\ol E_f$ be the support of the spectral density $f$ of a stationary process $X(t)$.
Then for the the minimal eigenvalue $\lambda_{1,n}(f)$ of $T_n(f)$ the following asymptotic estimates hold.
\begin{itemize}
\item[{\rm(a)}] If $\ol E_f= \G_{2\al}(\theta_0)$, where $\G_{2\al}(\theta_0)$ is as in Example \ref{ex1}, then
\beq
\label{aaa2}
\limsup_{n\to\f}\sqrt[n]{\lambda_{1,n}(f)} \leq \sin^{2}\left(\alpha/ 2\right).
\eeq
\item[{\rm(b)}] If $\ol E_f= \G_{\al,\de}(\theta_0)$, where $\G_{\al,\de}(\theta_0)$ is as in Example \ref{ex3}, then
\beq
\label{aaa3}
\limsup_{n\to\f}\sqrt[n]{\lambda_{1,n}(f)} \leq \sin(\al/2)\sin(\al/2+\de).
\eeq
\end{itemize}
\end{thm}
It is important to note that in Theorem \ref{T6.5} the essential
infimum $m_f=0$
is attained at uncountably many points, and the spectral density $f$,
in general, is not continuous on $E_f$ or at the endpoints of $E_f$.
Thus, $f$ does not satisfy Condition A and yet the rate of convergence
of $\lambda_{1,n}(f)$ to zero is much faster than in (\ref{ea3}).

Using Davisson's theorem (Theorem \ref{D1}) and its extension
(Theorem \ref{D1E}) we obtain exact upper bounds for the minimal
eigenvalue $\lambda_{1,n}(f)$ rather than the asymptotic estimates
\eqref{aaa2} and \eqref{aaa3}.

\begin{thm}[Babayan and Ginovyan \cite{BG-2}]
Let $\ol E_f$ be the support of the spectral density $f$ of a stationary process $X(t)$.
Then for the the minimal eigenvalue $\lambda_{1,n}(f)$ the following
inequalities hold.
\begin{itemize}
\item[{\rm(a)}] If $\ol E_f= \G_{2\al}(\theta_0)$, where $\G_{2\al}(\theta_0)$ is as in Example \ref{ex1}, then in view of \eqref{ndd1} and the first inequality in \eqref{ea7} we have
\beq
\label{dndd1}
\lambda_{1,n}(f) \le 4c\left(\sin(\alpha/2)\right)^{2n-2},
\eeq
where $c=r(0)$ and $r(\cdot)$ is the covariance function of $X(t)$.
\item[{\rm(b)}] If $\ol E_f= \G_{\al,\de}(\theta_0)$, where $\G_{\al,\de}(\theta_0)$ is as in Example \ref{ex3}, then in view of \eqref{d1e} and the first inequality in \eqref{ea7} we have
\beq
\label{dd1e}
\lambda_{1,n}(f)\le 4c\left(\sin(\alpha/2)\right)^{n-1}\left(\sin(\alpha/2+\de)\right)^{n-1}.
\eeq
\end{itemize}
\end{thm}
\ssn{A model with spectral density $f$ possessing exponential order zeros.}
Let $X(t)$ be a stationary process with spectral density $f_a$ given by
formula \eqref{nd4}, that is, $f_a$ is the Pollaczek-Szeg\H{o} function.
As it was observed (see \eqref{tt2})
\beq \label{tt222}
\nonumber
f_a(\la)\sim
\left \{
\begin{array}{ll}
 2e^a\exp\left\{-{a\pi}/{|\la|}\right\} & \mbox{as $\la\to0$},\\
2\exp\left\{-{a\pi}/{(\pi-|\la|)}\right\} & \mbox{as $\la\to\pm\pi$}.
\end{array}
\right.
\eeq
Thus, the function $f_a$ in \eqref{tt222} has a zero at points
$x=0,\pm\pi$ of exponential order and is positive elsewhere ($m_f=0<M_f$).
Then by Theorem \ref{R2} we have
\bea
\label{ea14}
\de_n(f_a) = \si^2_n(f_a)\sim n^{-a} \q {\rm as} \q n\to\f.
\eea
Now, by using the first inequality in \eqref{ea7}, from \eqref{ea14}
we conclude that
\bea
\label{ea15}
\la_{1,n}(f_a) = O\left(n^{-a}\right) \q {\rm as} \q n\to\f.
\eea
Thus, by choosing $a$ large enough one can obtain very fast rate of convergence
of $\la_{1,n}(f_a)$ to zero as $n\to\f$.

\begin{rem}
{\rm The asymptotic relation \eqref{ea15} remains valid for more
general models. Indeed, let $X(t)$ be a stationary process with
spectral density given by $f(\la)=f_a(\la)g(\la)$,
where $f_a(\la)$ is as in \eqref{nd4}, and $g(\la)$ is a function
satisfying the conditions of one of Theorems
\ref{sT1}, \ref{T1}, \ref{T2}, \ref{T1b} or Corollary \ref{ec5.2}.
Then for the
minimal eigenvalue $\lambda_{1,n}(f)$ of a truncated Toeplitz matrix
$T_n(f)$ generated by spectral density $f$, we have
\bea
\label{ea115}
\nonumber
\la_{1,n}(f) = O\left(n^{-a}\right) \q {\rm as} \q n\to\f.
\eea}
\end{rem}
More results concerning asymptotic behavior of the extreme
eigenvalues of truncated Toeplitz matrices can be found in
B\"ottcher and Grudsky \cite{BoG}, and Serra \cite{Se1}--\cite{Se3}).

Observe that in the above models information from the theory of
stationary processes was used to find linear-algebra results.
It is of interest to study the inverse problem.

\section{Tools}
\label{methods}

In this section we briefly discuss the tools, used to prove
the results stated in Sections \ref{ndp} and \ref{d}
(see Babayan and Ginovyan \cite{BG-1} -- \cite{BG-2}, and
Babayan et al. \cite{BGT}.

\sn{Some metric characteristics of bounded closed sets in the plane}
\label{DCC}
We introduce here some metric characteristics of bounded closed sets in the plane, such as, the transfinite diameter, the Chebyshev constant,
and the capacity, and discuss some properties of these characteristics.  Then, we state the theorems of Fekete and Robinson on the transfinite diameters of related sets, as well as, an extension of Robinson's theorem.

\ssn{The transfinite diameter, the Chebyshev constant and the capacity}

One of the fundamental result of geometric complex analysis is the
classical theorem by Fekete and Szeg\H{o}, stating that for any compact set $F$
in the complex plane  $\mathbb{C}$ the transfinite diameter, the Chebyshev constant
and the capacity of $F$ {\em coincide}, although they are defined from very
different points of view. Namely, the transfinite diameter of the set $F$ characterizes
the asymptotic size of $F$, the Chebyshev constant of $F$ characterizes the
minimal uniform deviation of a monic polynomial on $F$, and the capacity of $F$
describes the asymptotic behavior of the Green function at infinity.
For the definitions and results stated in this subsection we refer the
reader to the following references: Fekete \cite{F}, Goluzin \cite{GMG},
Chapter 7, Kirsch \cite{Kir}, Landkof \cite{La}, Chapter II, Ransford
\cite{Ran}, Chapter 5, Saff \cite{Saf}, Szeg\H{o} \cite{Sz1}, Chapter 16, and Tsuji \cite{Tsu}, Chapter III.

\vskip2mm
\n \underline{{\it Transfinite diameter.}}
Let $F$ be an infinite bounded closed (compact) set in the complex plane $\mathbb{C}$.
Given a natural number $n\geq 2$ and points $z_1, \ldots, z_n \in F$, we define
\beq
\label{td1}
d_n(F):=\max_{z_1,\ldots,z_n\in F}\left[\prod_{1\le j<k
\le n}^n|z_j-z_k|\right]^{2/[n(n-1)]},
\eeq
which is the maximum of products of distances between the
$ \begin{pmatrix} n\\2\end{pmatrix}=n(n-1)/2$
pairs of points $z_k$, $k=1,\ldots,n$, as the points $z_k$ range over the set $F$.
Note that $d_2(F)$ is the diameter of $F$.
The quantity $d_n(F)$ is called the $n$th
{\it transfinite diameter} of the set $F$.
It can be shown (see, e.g., Goluzin \cite{GMG}, p. 294)
that the sequence $d_n(F)$ is non-increasing and does not exceed the
diameter $d_2(F)$ of $F$, implying that $d_n(F)$ has a finite limit as $n\to\f$.
This limit, denoted by $d_\f(F)$, is called the {\it transfinite diameter} of
$F$. Thus, we have
\beq
\label{td2}
d_\f(F):=\lim_{n\to\f}d_n(F).
\eeq
If $F$ is empty or consists of a finite number of points, we put $d_\f(F)=0$.

\vskip2mm
\n \underline{{\it Chebyshev constant.}}
Let $F$ be as before, we put $m_n(F): = \inf\max_{z\in F}|q_n(z)|$,
where the infimum is taken over all monic polynomials $q_n(z)$ from the
class $\mathcal{Q}_n$, where $\mathcal{Q}_n$ is as in \eqref{Q_n}.
Then there exists a unique monic polynomial $T_n(z,F)$ from the class $\mathcal{Q}_n$,
called the {\it Chebyshev polynomial} of $F$ of order $n$, such that
\beq
\label{td6}
m_n(F)=\max_{z\in F}|T_n(z,F)|.
\eeq
Fekete \cite{F} proved that $\lim_{n\to\f}(m_n(F))^{1/n}$ exists.
This limit, denoted by $\tau(F)$, is called the {\it Chebyshev constant}
for the set $F$. Thus,
\beq
\label{td7}
\tau(F): = \lim_{n\to\f}(m_n(F))^{1/n}.
\eeq

\vskip2mm
\n \underline{{\it Capacity (logarithmic)}}.
Let $F$ be as above,
and let $D_F$ denote the complementary domain to $F$, containing
the point $z=\f$.
If the boundary $\G:=\partial D_F$ of the domain $D_F$
consists of a finite number of rectifiable Jordan curves, then for the domain $D_F$
can be constructed a Green function $G_{F}(z,\f):=G_{D_F}(z,\f)$ 
with a pole at infinity.
This function is harmonic everywhere in $D_F$, except at the point $z=\f$,
is continuous including the boundary $\G$ and vanishes on $\G$.
It is known that in a vicinity of the point $z=\f$ the function $G_{F}(z,\f)$ admits
the representation (see, e.g.,  Goluzin \cite{GMG}),  p. 309-310):
\beq
\label{tdG10}
G_F(z,\f)= \ln|z|+\g+ O(z^{-1}) \q {\rm as} \q z\to\f.
\eeq
The number $\g$ in \eqref{tdG10} is called the {\it Robin's constant} of the domain $D_F$,
and the number
\beq
\label{tdG11}
C(F):= e^{-\g}
\eeq
is called the  {\it capacity} (or the {\it logarithmic capacity}) of the set $F$.

\vskip2mm
Now we are in position to state the above mentioned fundamental result of
geometric complex analysis, due to M. Fekete and G. Szeg\H{o}
(see, e.g., Goluzin \cite{GMG}, p. 197 
and Tsuji \cite{Tsu}, p. 73).

\begin{pp}[Fekete - Szeg\H{o}'s theorem]
\label{FS1}
For any compact set $F\subset \mathbb{C}$, the transfinite diameter $d_\f(F)$
defined by \eqref{td2}, the Chebyshev constant $\tau(F)$ defined by \eqref{td7},
and the capacity $C(F)$ defined by \eqref{tdG11} coincide, that is,
\beq
\label{td11}
d_\f(F)=C(F)=\tau(F).
\eeq
\end{pp}
It what follows, we will use the term "transfinite diameter" and the notation $\tau(F)$
for \eqref{td11}.

The calculation of the transfinite diameter
(and hence, the capacity and the Chebyshev constant) is a challenging problem, and in only very few cases has the transfinite diameter been exactly calculated (see, e.g., Landkof \cite{La}, p. 172-173, Ransford \cite{Ran}, p.135, and also Examples \ref{ex1} -- \ref{ex4}).

In the next proposition we list a number of properties of the transfinite
diameter (and hence, of the capacity and the Chebyshev constant),
which were used to prove the results stated in Section \ref{d1}.
\begin{pp}
\label{pp1}
The transfinite diameter possesses the following properties.
\begin{itemize}
\item[(a)]
The transfinite diameter is monotone, that is, for any closed sets $F_1$ and $F_2$ with
$F_1\subset F_2$, we have $\tau(F_1)\leq \tau(F_2)$ (see, e.g., Saff \cite{Saf}, p. 169,
Tsuji \cite{Tsu}, p. 56).
\item[(b)]
If a set $F_1$ is obtained from a compact set $F\subset \mathbb{C}$ by
a linear transformation, that is, $F_1:=aF+b=\{az+b :\, z\in F\}$, then
$\tau(F_1)=|a|\tau(F)$. In particular, the transfinite diameter $\tau(F)$
is invariant with respect to parallel translation and rotation of $F$
(see, e.g., Goluzin \cite{GMG}, p. 298, Saff \cite{Saf}, p. 169,
Tsuji \cite{Tsu}, p. 56).
\item[(c)]
The transfinite diameter of an arbitrary circle of radius $R$
is equal to its radius $R$. In particular, the transfinite diameter of the unit
circle $\mathbb{T}$ is equal to 1 (Tsuji \cite{Tsu}, p. 84).
\item[(d)]
The transfinite diameter of an arc $\G_\al$ of a circle of radius $R$
with central angle $\al$ is equal to $R\sin(\al/4)$.
In particular, for the unit circle $\mathbb{T}$, we have
$\tau(\G_\al)=\sin(\al/4)$ (Tsuji \cite{Tsu}, p. 84).
\item[(e)]
The transfinite diameter of an arbitrary line segment $F$ is equal
to one-fourth its length, that is, if $F:=[a, b]$, then
$\tau(F) = \tau([a, b])=(b-a)/4.$ (see, e.g., Tsuji \cite{Tsu}, p. 84).
\end{itemize}
\end{pp}

\n \underline{{\it The inner and outer transfinite diameters. $\tau$-measurable sets}}.
For an arbitrary (not necessarily closed) bounded set $E\subset \mathbb{C}$, we also define the {\em inner} and the {\em outer transfinite diameters}, denoted by $\tau_*(E)$ and $\tau^*(E)$, respectively,
as follows (see, e.g., Babayan \cite{Bb-1}, Korovkin \cite{Kor1}):
\beq
\label{tdo1}
\tau_*(E): = \sup_{F\subset E}\tau(F) \quad {\rm and} \quad
\tau^*(E): = \tau(\ol E),
\eeq
where the supremum is taken over all compact subsets $F$ of the set $E$,
and $\ol E$ stands for the closure of $E$.
Observe that $\tau_*(E) \le \tau^*(E)$. The set $E$ for which
$\tau_*(E) = \tau^*(E)$ is said to be $\tau$-{\em measurable},
and in this case, we write $\tau(E)$ for the common value: $\tau(E)=\tau_*(E) = \tau^*(E).$

\ssn{Transfinite diameters of related sets}
\label{FR}

We state here the theorems of Fekete \cite{F} and  Robinson \cite{Rob}
on the transfinite diameters of related sets, as well as,
an extension of  Robinson's theorem.

The following classical theorem about the transfinite diameter
of related sets was proved by Fekete \cite{F} (see also Goluzin
\cite{GMG}, pp. 299-300).

\begin{pp}[Fekete's theorem, Fekete \cite{F}]
\label{FT}
Let $F$ be a bounded closed set in the complex $w$-plane $\mathbb{C}_w$,
and let $p(z):=p_n(z)=z^n+c_1z^{n-1}\cdots+c_n$ be an arbitrary monic
polynomial of degree $n$. Let $F^*$ be the preimage of $F$
in the $z$-plane $\mathbb{C}_z$ under the mapping $w=p(z)$, that is,
$F^*$ is the set of all points $z\in\mathbb{C}_z$ such that $w:=p(z)\in F$.
Then
\beq
\label{F0}
\tau(F^*)=[\tau(F)]^{1/n},
\eeq
where $\tau(F)$ and $\tau(F^*)$ stand for the transfinite diameters
of the sets $F$ and $F^*$, respectively.
\end{pp} 
Observe that if in Fekete's theorem the mapping is carried out by
an arbitrary (not necessarily monic) polynomial of degree $n$:
$p(z):=p_n(z)=az^n+a_1z^{n-1}\cdots+a_n$ ($a\neq0$), then we have
\beq
\label{F00}
\tau(F^*)=[\tau(F)/|a|]^{1/n}.
\eeq
In Robinson \cite{Rob}, Fekete's theorem was extended
to the case where the mapping is carried out by a rational function
instead of a polynomial. More precisely, in Robinson \cite{Rob}
was proved the following theorem.

\begin{pp}[Robinson's theorem, Robinson \cite{Rob}]
\label{RT1}
Let $p(z):=p_n(z)=z^n+a_1z^{n-1}\cdots+a_n$ and $q(z):=q_k(z)$ be arbitrary
relatively prime polynomials of degrees $n$ and $k$, respectively,
with $k<n$.
Let $F$ be a bounded closed set in the complex $w$-plane $\mathbb{C}_w$,
and let $F^*$ be the preimage of $F$ in the $z$-plane $\mathbb{C}_z$
under the mapping $w=\I(z):=p(z)/q(z)$. Assume that $|q(z)|=1$ for all
$z\in F^*$. Then
\beq
\label{F9}
\tau(F^*)=[\tau(F)]^{1/n}.
\eeq
\end{pp}
\begin{rem}
\label{mk1}
{\rm It is clear that the condition $|q(z)|=1$ for all $z\in F^*$
in Robinson's theorem can be
replaced by the condition $|q(z)|=C$ for all $z\in F^*$ with
an arbitrary positive constant $C$, and, in this case, the
relation \eqref{F9} becomes $\tau(F^*)=[C\tau(F)]^{1/n}$.}
\end{rem}

Observe that in the special case where $q(z)\equiv 1$, Robinson's
theorem reduces to the Fekete theorem (Proposition \ref{FT}).
A special interest represents the other special case where
$p(z)=z^2+1$ and $q(z)=2z$. In this case, the mapping given by the
rational function
$$\I(z):=\frac{p(z)}{q(z)}=\frac12\left(z+\frac1z\right)$$
projects the subsets of the unit circle $\mathbb{T}:=\{z\in\mathbb{C}:\, |z|=1\}$
onto the real axis $\mathbb{R}$, and, in view of Remark \ref{mk1},
Robinson's theorem (Proposition \ref{RT1}) reads as follows.

\begin{pp}[Robinson \cite{Rob}]
\label{RT2}
Let $F$ be a bounded closed subset of the complex plane $\mathbb{C}$ lying
on the unit circle $\mathbb{T}$ and symmetric with respect to real axis,
and let $F^x$ be the projection of $F$ onto the real axis. Then
\beq
\label{F10}
\tau(F)=[2\tau(F^x)]^{1/2}.
\eeq
\end{pp}

\begin{rem}
\label{ExR}
{\rm The examples given in Section \ref{EX} show that the formula \eqref{F10} gives a
simple way to calculate the transfinite diameters of some
subsets of the unit circle $\mathbb{T}$, based only on the formula
of the transfinite diameter of a line segment (see Proposition \ref{pp1}(e)).}
\end{rem}

\ssn{An extension of Robinson's theorem}
\label{ER}

As it was observed above, the condition $|q(z)|=C$ for all $z\in F^*$
in Robinson's theorem (Proposition \ref{RT1}) is too restrictive,
and it essentially reduces the range of applicability of the theorem
into the following two cases:

(a) $q(z)\equiv 1$, and Robinson's theorem reduces to the Fekete
theorem.

(b) $p(z)=z^2+1$ and $q(z)=2z$, and, in this case,
the rational function $\I(z)=(z^2+1)/(2z)$
projects the subsets of the unit circle $\mathbb{T}$
onto the real axis $\mathbb{R}$.

Therefore, the question of extending Robinson's theorem
to the case where the condition $|q(z)|=M$
is replaced by a weaker condition becomes topical.
The next result, which was proved in Babayan and Ginovyan \cite{BG-2},
provides such an extension.
\begin{pp}[Babayan and Ginovyan \cite{BG-2}]
\label{MT}
Let the polynomials $p(z)$, $q(z)$, the sets $F$, $F^*$, and the mapping
$w=\I(z):=p(z)/q(z)$ be as in Proposition \ref{RT1}, and let $m:=\min_{z\in F^*}|q(z)|$ and
$M:=\max_{z\in F^*}|q(z)|$. Then the following inequalities hold:
\beq
\label{F11}
[m\tau(F)]^{1/n}\leq \tau(F^*)\leq [M\tau(F)]^{1/n}.
\eeq
\end{pp}
\begin{rem}
{\rm If the condition $|q(z)|=C$ is satisfied for all $z\in F^*$,
then we have $m=M =C$, and Proposition \ref{MT} reduces to Robinson's
theorem (Proposition \ref{RT1}).}
\end{rem}
\begin{rem}
{\rm Proposition \ref{MT} can easily be extended to more general case
where $p(z):=p_n(z)$ is an arbitrary (not necessarily monic) polynomial
of degree $n$:
$p(z)=az^n+a_1z^{n-1}\cdots+a_n$, $a\neq0$. Indeed, in this case,
canceling the fraction $\I(z):=p(z)/q(z)$ by $a$, we get
$\I(z):=p_1(z)/q_1(z)$, where now
$p_1(z)=p(z)/a=z^n + {\rm lower \,\,  order \,\, terms}$,
is a monic polynomial. Also, we have $\min_{z\in F^*}|q_1(z)|=m/|a|$
and $\max_{z\in F^*}|q_1(z)|=M/|a|$, where $m$ and $M$ are as in Proposition \ref{MT}.
Hence we can apply the inequality \eqref{F11} to obtain
\beq
\label{F16}
\nonumber
\left[\frac{m}{|a|}\tau(F)\right]^{1/n}\leq \tau(F^*)\leq \left[\frac{M}{|a|}\tau(F)\right]^{1/n}.
\eeq
}
\end{rem}

\sn{Weakly varying sequences}
\label{ww}
We introduce here the notion of weakly varying sequences and state
some of their properties (see Babayan et al. \cite{BGT}).
This notion was used in the specification of the class $\mathcal{F}$
of deterministic
processes considered in Section \ref{d2} (see formula \eqref{k04}).
\begin{den}
\label{kd1}
A sequence of non-zero numbers $\{a_n, \, n \in\mathbb{N}\}$ is said 
to be weakly varying if 
$$\lim_{n\to\f} {a_{n+1}}/{a_{n}} =1.$$
\end{den}
For example, the sequence $\{n^{\al}, \, \, \al\in\mathbb{R}, \, n\in\mathbb{N}\}$
is weakly varying (for $\al<0$ it is weakly decreasing and for $\al>0$
it is weakly increasing), while the geometric progression $\{q^n, \, 0<q<1, \,
n \in\mathbb{N}\}$ is not weakly varying.

In the next proposition we list some simple properties of the weakly varying sequences, which can easily be verified.
\begin{pp}
\label{p4.1}
The following assertions hold.
\begin{itemize}
\item[(a)]
If $a_n$ is a weakly varying sequence, then
$\lim_{n\to\f} {a_{n+\nu}}/{a_{n}} =1$ for any $\nu\in \mathbb{N}$.
\item[(b)]
If $a_n$ is such that $\lim_{n\to\f}a_n=a\neq0$,
then $a_n$ is a weakly varying sequence.
\item[(c)]
If $a_n$ and $b_n$ are weakly varying sequences, then $ca_n$ $(c\neq0),$
$a_n^\al \, (\al\in\mathbb{R}, a_n>0)$, $a_nb_n$ and $a_n/b_n$ also are
weakly varying sequences.
\item[(d)]
If $a_n$ is a weakly varying sequence,
and $b_n$ is a sequence of non-zero numbers such that
$\lim_{n\to\infty}{b_n}/{a_n}=c\neq 0,$
then $b_n$ is also a weakly varying sequence.
\item[(e)]
If $a_n$ is a weakly varying sequence of
positive numbers, then it is exponentially neutral
(see Definition \ref{ekd1}(a) and Remark \ref{vvv}).
\end{itemize}
\end{pp}

\sn{Some properties of the geometric mean and trigonometric polynomials.}
Recall that a trigonometric polynomial $t(\lambda)$ of degree $\nu$ is a
function of the form:
\beq
\label{tp1}
\nonumber
t(\lambda)= a_0+\sum_{k=1}^\nu(a_k\cos k\la + b_k\sin k\la)
= \sum_{k=-\nu}^\nu c_ke^{ik\la}, \q \la\in \mathbb{R},
\eeq
where $a_0, a_k, b_k \in \mathbb{R}$, $c_0=a_0$, $c_k=1/2(a_k-ib_k)$,
$c_{-k}=\ol c_k=1/2(a_k+ib_k)$, $k=1,2,\ldots, \nu$.

Recall that for a function $h\geq0$ by $G(h)$ we denote the geometric
mean of $h$ (see formula \eqref{a2}).
In the next proposition we list some properties of the geometric mean $G(h)$
and trigonometric polynomials (see Babayan et al. \cite{BGT}).
\begin{pp}
\label{p4.2}
The following assertions hold.
\begin{itemize}
\item[(a)]
Let $c>0$, $\al\in\mathbb{R}$, $f\geq0$ and $g\geq0$. Then
\beq
\label{gt1}
G(c)=c; \q G(fg)=G(f)G(g); \q G(f^\al)=G^\al(f) \, \, (G(f)> 0).
\eeq
\item[(b)]
$G(f)$ is a non-decreasing functional of $f$: if \ $0\leq f(\la)\leq g(\la)$,
then $0\leq G(f)\leq G(g)$.
In particular, if $0\leq f(\la)\leq 1$, then $0\leq G(f)\leq 1$.
\item[(c)]
{\rm (Fej\'er-Riesz theorem)}. Let $t(\la)$ be a non-negative trigonometric polynomial of degree $\nu$. Then
there exists an algebraic polynomial  $s_{\nu}(z)$ $(z\in \mathbb{C})$
of the same degree $\nu$, such that $s_{\nu}(z)\neq0$ for $|z|<1$, and
\beq
\label{ss2}
t(\lambda)=|s_{\nu}(e^{i\lambda})|^2.
\eeq
Under the additional condition $s_{\nu}(0)>0$ the polynomial $s_{\nu}(z)$
is determined uniquely.
\item[(d)]
Let $t(\lambda)$ and  $s_{\nu}(z)$ be as in Assertion (c). Then
$G(t)=|s_{\nu}(0)|^2>0$.
\end{itemize}
\end{pp}

\n
{\bf Acknowledgments.}

\n
The authors are grateful to their supervisor Academician, 
Professor I.A. Ibragimov
for introducing them to this research area.

\end{document}